\documentclass[11pt]{amsart}

\usepackage{amscd}
\usepackage{amssymb}

\usepackage{mathrsfs}

\swapnumbers
\usepackage{amsthm}

\newtheorem{thm}[equation]{Theorem}
\newtheorem{lem}[equation]{Lemma}
\newtheorem{cor}[equation]{Corollary}
\newtheorem{prop}[equation]{Proposition}

\newtheorem*{thm*}{Theorem}
\newtheorem*{prop*}{Proposition}
\newtheorem*{cor*}{Corollary}
\newtheorem*{lem*}{Lemma}

\theoremstyle{definition} %
\newtheorem{defn}[equation]{Definition}
\newtheorem*{defn*}{Definition}

\newtheorem{eg}[equation]{Example}
\newtheorem{egs}[equation]{Examples}

\theoremstyle{remark} %
\newtheorem{rmk}[equation]{Remark}

\newtheorem*{rmk*}{Remark}
\newtheorem*{rmks*}{Remarks}

\newtheoremstyle{exercise}
  {3pt}
  {3pt}
  {\small}
  {\parindent}
  {\sc\small}
  {.}
  {.5em}
   {}     
  {}

\theoremstyle{exercise}

\numberwithin{equation}{section}

\makeatletter
{\renewcommand{\theequation}{#1}}%
{\renewcommand{\theequation}{\arabic{equation}}\addtocounter{equation}{-1}\global\@ignoretrue}
\makeatother

\makeatletter
{\renewcommand{\theequation}{#1}\begin{eqnarray}}%
{\end{eqnarray}\renewcommand{\theequation}{\arabic{equation}}\addtocounter{equation}{-1}\global\@ignoretrue}
\makeatother

\makeatletter
\newenvironment{borel}[1]%
{\smallskip \refstepcounter{equation}\noindent{\textbf \theequation. }{{\textbf{#1.}}}}%
{\smallskip \global\@ignoretrue}
\makeatother

\makeatletter
{\smallskip \refstepcounter{equation}{\sc \theequation}{\sc (#1).}}%
{\smallskip \global\@ignoretrue}
\makeatother

\makeatletter
{\smallskip \refstepcounter{equation}\noindent{\sc \theequation.}{\sl{ #1.}}}%
{\smallskip \global\@ignoretrue}
\makeatother

\makeatletter
\newenvironment{borel*}%
{\smallskip \refstepcounter{equation}\noindent{\textbf \theequation.}}%
{\global\@ignoretrue}
\makeatother

\newcommand{\flist}[1]{\hangindent\leftmargini\textup{(1)}\hskip\labelsep {#1}%
\begin{enumerate}%
\setcounter{enumi}{1}%
}

\newcommand{\ot}{\otimes}

\newcommand{\C}{{\mathbb{C}}}        
\newcommand{\Q}{{\mathbb{Q}}}        
\newcommand{\R}{{\mathbb{R}}}        
\newcommand{\Z}{{\mathbb{Z}}}        
\newcommand{\N}{{\mathbb{N}}}      
\newcommand{\hyp}{\mathcal{H}}

\newcommand{\Zm}[1]{\Z/{#1}\Z}

\newcommand{\la}{\lambda}

\newcommand{\G}{{\Gamma}}       

\newcommand{\oddots}{{\mathinner{\mkern1mu\raise1pt\vbox{\kern7pt\hbox{.}}\mkern2mu\raise4pt\hbox{.}\mkern2mu\raise7pt\hbox{.}\mkern1mu}}}

\newcommand{\ksep}{k_{{\mathrm{sep}}}}

\newcommand{\Kx}{{K^{\times}}}

\newcommand{\kx}{k^\times}

\newcommand{\qform}[1]{{\langle{#1}\rangle}}                   
\newcommand{\pform}[1]{{\langle\!\langle{#1}\rangle\!\rangle}} 

\DeclareMathOperator{\Spin}{Spin}           

\newcommand{\oE}{\ensuremath{^1\!E_6}}
\newcommand{\dE}{\ensuremath{^2\!E_6}}

\DeclareMathOperator{\rank}{rank}

\DeclareMathOperator{\chr}{char}

\DeclareMathOperator{\Id}{Id}

\DeclareMathOperator{\cd}{cd}

\DeclareMathOperator{\End}{End}

\DeclareMathOperator{\aut}{Aut}

\newcommand{\Hom}{{\mathrm{Hom}}}

\newcommand{\iso}{\xrightarrow{\sim}}

\newcommand{\ra}{\rightarrow}

\theoremstyle{plain}

\renewcommand{\theenumi}{\alph{enumi}}

\newcommand{\iV}{{V^\iota}}

\DeclareMathOperator{\Str}{Str}

\DeclareMathOperator{\GL}{GL}
\DeclareMathOperator{\SL}{SL}
\DeclareMathOperator{\PSL}{PSL}

\newcommand{\br}[1]{[#1]}

\newcommand{\hrm}{\mathscr{H}}

\newcommand{\T}{\mathscr{T}}
\newcommand{\Td}{\T^d}

\newcommand{\ach}{\check{\alpha}}

\newcommand{\Ad}{A^d}
\newcommand{\Jd}{\Ad}

\theoremstyle{remark}

\newcommand{\g}{\mathfrak{g}}

\newcommand{\lradG}{0.25}
\newcommand{\darkradE}{0.115}

\begin{document}

\title%
[Groups of outer type $E_6$]{Groups of outer type $E_6$ with trivial Tits algebras}
\author{Skip Garibaldi}
\address{Department of Mathematics \& Computer Science, Emory University, Atlanta, GA 30322, USA}
\email{skip@member.ams.org}
\urladdr{http://www.mathcs.emory.edu/{\textasciitilde}skip/}

\author{Holger P. Petersson}
\address{Fachbereich Mathematik, FernUniversit\"at in Hagen, D-58084 Hagen, Germany}
\email{Holger.Petersson@FernUni-Hagen.de}
 \urladdr{http://www.fernuni-hagen.de/MATHEMATIK/ALGGEO/Petersson/petersson.html}

\begin{abstract}
In two 1966 papers, J.~Tits gave a construction of exceptional Lie algebras (hence implicitly exceptional algebraic groups) and a classification of possible indexes of simple algebraic groups.  For the special case of his construction that gives groups of type $E_6$, we connect the two papers by answering the question: Given an Albert algebra $A$ and a separable quadratic field extension $K$, what is the index of the resulting algebraic group?
\end{abstract}

\date{\today}

\maketitle

\setcounter{tocdepth}{1}
\tableofcontents

\newsavebox{\dEpic}
\savebox{\dEpic}(5,1){\begin{picture}(5, 1)
    \multiput(3,0.75)(1,0){2}{\circle*{\darkradE}}
    \multiput(3,0.25)(1,0){2}{\circle*{\darkradE}}
    \multiput(1,0.5)(1,0){2}{\circle*{\darkradE}}
    
    \put(1, 0.5){\line(1,0){1}}
    \put(3,0.75){\line(1,0){1}}
    \put(3, 0.25){\line(1,0){1}}
    
    \put(3,0.5){\oval(2,0.5)[l]}\end{picture}}
    
The exceptional simple algebraic groups are organized in a chain of inclusions
\[
A_1 \subset A_2 \subset G_2 \subset D_4 \subset F_4 \subset E_6 \subset E_7 \subset E_8.
\]
One approach to proving something about a group of type, say, $E_6$, is to attempt to make use of known facts about the groups of types appearing earlier in the chain.  Essentially everything is known about groups of types $A_1$ (corresponding to quaternion algebras), $A_2$ 
by \cite[\S19]{KMRT}, and $G_2$ (corresponding to octonion algebras).  Quite a lot is known about groups of type $F_4$, corresponding to Albert algebras.  In contrast, very little is known about groups of type $E_6$.  Two of the main results are Tits's construction in \cite{Ti:const} and the classification of possible indexes \cite{Ti:Cl}.

The version of Tits's construction studied here takes an Albert algebra $A$ and a quadratic \'etale algebra $K$ and produces a simply connected group $G(A, K)$ of type $E_6$.   (See \S\ref{remind} below for background on Albert algebras.) The purpose of this paper is to answer the question:
\begin{equation} \label{ques}
\text{\emph{What is the index of the group $G(A, K)$?}}
\end{equation}
The hard part of answering this question is treated by the following theorem.  We fix an arbitrary base field $k$.  

\begin{thm}  \label{MT}
The following are equivalent:
\begin{enumerate}
\renewcommand{\theenumi}{\arabic{enumi}}
\item The group $G(A, K)$ is isotropic.
\item $k \times K$ is (isomorphic to) a subalgebra of $A$.
\item $A$ is reduced and there exists a $2$-Pfister bilinear form $\gamma$ such that 
\[ \gamma \cdot f_3(A) = f_5(A) \quad \text{and} \quad  \gamma \cdot [K] = 0.
\]
\end{enumerate}
Conditions $(1)$--$(3)$ are implied by:
\begin{enumerate}
\item[(4)] $A$ has a nonzero nilpotent element.
\end{enumerate}
Furthermore, if $A$ is split by $K$, then $(1)$--$(3)$ are equivalent to $(4)$.
\end{thm}

Once one knows that $G(A, K)$ is isotropic, it is not difficult to determine the index of $G(A, K)$---see Prop.~\ref{trivial.prop} and \ref{rGAK}---so we have completely settled Question \eqref{ques}.

When $K$ is ``split'' (i.e., $K = k \times k$, equivalently, $G(A, K)$ has type $\oE$), the theorem is a triviality.  Indeed, conditions (1) through (3) are equivalent to the statement ``$A$ is reduced''.  When $K$ is split, we  define the statement ``$A$ is split by $K$'' to mean that $A$ is split as an Albert $k$-algebra.

The main theorem shows the flavor of the paper; it mixes algebraic groups (in (1)), Jordan algebras (in (2) and (4)), and---essentially---quadratic forms (in (3)).   The core of our proof is Jordan-theoretic.  We prove 
(1) implies (2) or (4) in Cor.~\ref{transl} and Propositions \ref{Tx0.prop} and \ref{Tx1.prop}.  We prove that (4) implies (2) in \ref{imp42}, (2) implies (3) in \ref{imp23}, and (3) implies (1) in \ref{iso.cor}.  The last claim  is proved in Example \ref{nil.eg2} and \ref{furthermore}.

As side benefits of the proof, we obtain concrete descriptions of the projective homogeneous varieties for groups of type $\dE$ in \S\ref{homo} and we easily settle an open question from a 1969 paper of Veldkamp \cite{Veld:unit2} in \ref{Veld}.

\section{Notation and reminders} 

Recall that the \emph{(Tits) index} \cite[2.3]{Ti:Cl} of an (affine, semisimple) algebraic group is its Dynkin diagram plus two other pieces of information: the Galois action on the diagram and circles indicating the maximal $k$-split torus in the group.  

The \emph{Tits algebras} of an algebraic group $G$ are the $k$-algebras $\End_G(V)$ as $V$ varies over $k$-irreducible representations of $G$.  We say that $G$  \emph{has trivial Tits algebras} if $\End_G(V)$ is a (commutative) field for every $V$.

Below, $K$ will always denote a quadratic \'etale $k$-algebra with nontrivial $k$-automorphism $\iota$, and $A$ is an Albert $k$-algebra.

Throughout, we use the notation $\qform{\alpha_1, \ldots, \alpha_n}$ for the diagonal matrix with $\alpha_i$ in the $(i, i)$ entry and for the symmetric bilinear form with that Gram matrix.  We write $\pform{\alpha_1, \ldots, \alpha_n}$ for the Pfister bilinear from $\qform{1, -\alpha_1} \ot \cdots \ot \qform{1, -\alpha_n}$.  \emph{Pfister form} means ``Pfister quadratic form".  We write $I^n k$ for the module of quadratic forms generated by the $n$-Pfister forms over the Witt ring of symmetric bilinear forms (this agrees with the usual notation in characteristic $\ne 2$).  We write $[K]$ for the 1-Pfister form given by the norm $K \ra k$;  a similar convention applies to the norm of an octonion $k$-algebra.

For $n \in \N$, we write $H^q(k, n)$ for the group denoted by $H^q(k, \Zm{n}(q - 1))$ in \cite[pp.~151--155]{GMS}.  When $n$ is not divisible by the characteristic of $k$, it is the Galois cohomology group $H^q(k, \boldsymbol{\mu}_n^{\ot (q - 1)})$.  There is a bijection between $q$-fold Pfister forms (up to isomorphism) and symbols in $H^q(k, 2)$.  In characteristic different from 2, this is a direct consequence of Voevodsky's proof of the Milnor Conjecture, and in characteristic 2 it is in \cite{ABa}.  We will write, for example, $[K]$ also for the symbol in $H^1(k, 2)$ corresponding to the norm $K \ra k$, and similarly for an octonion $k$-algebra.

\section{Rost invariants} \label{iso}

Let $G$ be a quasi-simple, simply connected group over a field $k$.  There is a canonical map $r_G \!: H^1(k, G) \ra H^3(k, n_G)$ known as the \emph{Rost invariant}, where $n_G$ is a natural number depending on $G$, see \cite{GMS} for details.  The map is ``functorial in $k$''.  In this section, we relate the Rost invariant with the index of isotropic groups of type $\dE$ with trivial Tits algebras.

\smallskip

There is a class $\nu \in H^1(k, \aut(G)^\circ)$ such that the twisted group $G_\nu$ is quasi-split.  Moreover, this property uniquely determines $\nu$, as can be seen from a twisting argument and the fact that the kernel of the map $H^1(k, \aut(G_\nu)^\circ) \ra H^1(k, \aut(G_\nu))$ is zero.

If $G$ has trivial Tits algebras, then there is an $\eta \in H^1(k, G)$ that maps to $\nu$ and we 
 write $a(G) \in H^3(k, n_G)$ for the Rost invariant $r_G(\eta)$.
\begin{lem}
The element $a(G)$ depends only on the isomorphism class of $G$ (and not on the choice of $\eta$).
\end{lem}

\begin{proof}
Fix a particular $\eta$.  Every inverse image of $\nu$ is of the form $\zeta \cdot \eta$ for some $\zeta \in Z^1(k, Z(G))$.  Write $\tau$ for the ``twisting" isomorphism $H^1(k, G) \iso H^1(k, G_\eta)$.  The centers of $G$ and $G_\eta$ are canonically identified, and we have
\[
\tau(\zeta \cdot \eta) = \zeta \cdot \tau(\eta) = \zeta \cdot 1.
\]
Since the Rost invariant is compatible with twisting,
\[
r_G(\zeta \cdot \eta) = r_{G_\eta}(\tau(\zeta \cdot \eta)) + r_G(\eta) = r_{G_\eta}(\zeta \cdot 1) + r_G(\eta).
\]
But $G_\eta$ is quasi-split, so the image $\zeta \cdot 1$ of $\zeta$ in $H^1(k, G_\eta)$ is trivial.
\end{proof}

\begin{egs} \label{trivial.eg}
In the following examples, $G$ always denotes a quasi-simple, simply connected group with trivial Tits algebras. 

\begin{enumerate}
\item Let $G$ be of type $^1\!D_n$ for $n = 3$ or 4, so that $G$ is isomorphic to $\Spin(q)$ for some $2n$-dimensional quadratic form $q$ in $I^3 k$.  The invariant $a(G)$ is the Arason invariant of $q$.

In the case $n = 3$, the Arason-Pfister Hauptsatz implies that $q$ is hyperbolic, so $G$ is split and $a(G)$ is zero.

In the case $n = 4$, $q$ is similar to a 3-Pfister form.  It follows that $a(G)$ is zero if and only if $G$ is split if and only if $G$ is isotropic.

\item Let $G$ be of type $^2\!D_n$ for $n = 3$ or 4, with associated separable quadratic extension $K/k$.  
The group $G$ is isomorphic to $\Spin(q)$ for $q$ a $2n$-dimensional quadratic form such that $q - [K]$ is in $I^3 k$; the Arason invariant of $q - [K]$ is $a(G)$.

We claim that $a(G)$ is a symbol in $H^3(k, 2)$.  
When $n = 3$, $q - [K]$ is an 8-dimensional form in $I^3 k$, so it is similar to a Pfister form.
In the case $n = 4$, $q - [K]$ is a 10-dimensional form in $I^3 k$, hence it is isotropic \cite[4.4.1(ii)]{Ti:si}.  The Hauptsatz implies that $q - [K]$ is isomorphic to $\qform{\alpha} \gamma \perp \hyp$ for some $\alpha \in \kx$ and some 3-Pfister $\gamma$, where $\hyp$ denotes a hyperbolic plane. 
This shows that $a(G)$ is a symbol in both cases.

We next observe that $a(G)$ is \emph{not} killed by $K$ if and only if $n = 4$ and $G$ is $k$-anisotropic.  To see this, we may assume that $n = 4$ and $G$ is $k$-anisotropic by (a); we will show that $G$ is anisotropic over $K$.  Suppose not, i.e., that $G$ is split by $K$, that is, $q$ is hyperbolic over $K$.  It follows that $\gamma$ is isomorphic to $\beta [K]$ for some 2-Pfister bilinear form $\beta$ \cite[4.2(iii)]{HL}.  Since $\qform{\alpha} \gamma \perp [K]$ is isomorphic to $q \perp \hyp$, $\gamma$ represents $-\alpha \la$ for some nonzero $\la$ represented by $[K]$.  The roundness of $\gamma$ and $[K]$ gives that $\qform{\alpha} \gamma$ is isomorphic to $\qform{-1}  \beta [K]$.   Thus, $q \perp \hyp$ has Witt index at least 2.  This contradicts the hypothesis that $G$ is $k$-anisotropic, completing the proof of the observation.

\item Let $G$ be \emph{anisotropic} of type $^2\!A_5$ with associated separable quadratic extension $K/k$; $G$ is isomorphic to the special unitary group of a $K/k$-hermitian form deduced from a 6-dimensional symmetric bilinear form $\beta$ over $k$.   Note that $a(G)$ lives in $H^3(k, 2)$ and is the Arason invariant of $\beta [K]$.  

Suppose that $a(G)$ is a symbol.  Since $G$ is split by $K$, $a(G)$ is of the form $[K] \cdot (\la) \cdot (\mu)$ for some $\la, \mu \in \kx$; that is, $\beta [K]$ is congruent to the corresponding 3-Pfister $\gamma$ modulo $I^4 k$.  The function field of $\gamma$ makes the 12-dimensional form $\beta [K]$ hyperbolic by the Hauptsatz, so the anisotropic part of $\beta [K]$ is isomorphic to $\tau \ot \gamma$ for some bilinear form $\tau$.  In particular, $\beta [K]$ is isotropic, contradicting the anisotropy of $G$.  We conclude that $a(G)$ is \emph{not} a symbol.
\end{enumerate}
\end{egs}

The purpose of the preceding examples was to prepare the proof of the following proposition.

\begin{prop} \label{trivial.prop}
Let $G$ be an isotropic group of type $\dE$ with trivial Tits algebras; write $K$ for the associated quadratic extension of $k$.  Then $a(G)$ is in $H^3(k, 2)$ and the index of $G$ is given by Table \ref{iso.table}.
\end{prop}

\refstepcounter{equation} \label{iso.table}
\begin{table}[thb]
\begin{tabular}{cc} \\
index & condition \\ \hline
quasi-split & $a(G) = 0$ \\
\parbox{5cm}{\begin{picture}(5, 1.1)
\put(0,0){\usebox{\dEpic}} 
\put(4,0.5){\oval(0.4,0.75)}
\put(1,0.5){\circle{\lradG}}
\end{picture}}& $a(G)$ is a nonzero symbol killed by $K$\\
\parbox{5cm}{\begin{picture}(5, 1.1)
\put(0,0){\usebox{\dEpic}} 
\put(4,0.5){\oval(0.4,0.75)}
\end{picture}}& $a(G)$ is a symbol not killed by $K$\\
\parbox{5cm}{\begin{picture}(5, 1)
\put(0,0){\usebox{\dEpic}} 
\put(1,0.5){\circle{\lradG}}
\end{picture}}& $a(G)$ is not a symbol
\end{tabular}
\caption{Tits indexes and their corresponding Rost invariants}
\end{table}

\begin{proof}
Consulting the list of possible indexes from \cite[p.~59]{Ti:Cl}, the only one missing from the table is
   \[
\parbox{5cm}
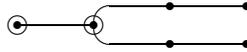
{\begin{picture}(5, 1)
\put(0,0){\usebox{\dEpic}} 
\put(1,0.5){\circle{\lradG}}
\put(2,0.5){\circle{\lradG}}
\end{picture}}
\]
But a group with such an index has nontrivial Tits algebras \cite[5.5.5]{Ti:R}, therefore the index of $G$ is in the table.

If $G$ is not quasi-split, then the semisimple anisotropic kernel $H$ of $G$ is quasi-simple with trivial Tits algebras, and it follows from Tits's Witt-type theorem that $a(G)$ equals $a(H)$.  The correspondence between the index and $a(G)$ asserted by the table now follows from Example \ref{trivial.eg}.
\end{proof}

A group $G$ from the first three rows of the table is completely determined by the value of $a(G)$.

\begin{prop} \label{class.prop}
Let $G$ and $G'$ be quasi-simple, simply connected groups of type $\dE$ whose indexes are in the first three rows of Table \ref{iso.table}.  If $a(G)$ equals $a(G')$, then $G$ and $G'$ are isomorphic.
\end{prop}

\begin{proof}
Fix a maximal $k$-split torus $S$ in $G$, a maximal $k$-torus $T$ containing it, and a set of simple roots for $G$ with respect to $T$.  Since $G$ is simply connected, the group of cocharacters of $T$ is identified with the coroot lattice.  Write $S_1$ for the rank 1 torus corresponding to $\ach_1 + \ach_6$ (``corresponding to the circle around the $\alpha_1$ and $\alpha_6$ vertices in the index"); it is $k$-defined by \cite[Cor.~6.9]{BoTi}.  Put $G_1$ for the derived subgroup of $Z_G(S_1)$; it is simply connected of type $^2\!D_4$ with trivial Tits algebras.

Define a subgroup $G'_1$ of $G'$ in an analogous manner.  Since $G_1$ and $G'_1$ are strongly inner forms of each other, 
$G'_1$ is isomorphic to $G_1$ twisted by some 1-cocycle $\alpha \in Z^1(k, G_1)$.  The semisimple anisotropic kernel of $G'$ lies in $G'_1$, hence is the semisimple anisotropic kernel of $G$ twisted by $\alpha$.  Tits's Witt-type theorem implies that $G'$ is isomorphic to $G$ twisted by $\alpha$.

Since $a(G)$ equals $a(G')$, the Rost invariant $r_G(\alpha)$ is zero by a twisting argument.  But the inclusion of $G_1$ into $G$ arises from the natural inclusion of root systems and so has Rost multiplier one.  Hence  $r_{G_1}(\alpha)$ is also zero.  Moreover, since $G_1$ has trivial Tits algebras, $r_{G_1}(\zeta \cdot \alpha)$ is zero for every 1-cocycle $\zeta$ with values in the center of $G_1$.

Fix an isomorphism of $G_1$ with $\Spin(q)$ for some 8-dimensional quadratic form $q$ and write $q_{\zeta \cdot \alpha}$ for the quadratic form obtained by twisting $q$ via the image of $\zeta \cdot \alpha$ in $Z^1(k, SO(q))$; the forms so obtained are precisely the forms $\qform{\la} q$ for $\la \in \kx$.  Pick a $\zeta$ (and hence a $\la$) so that $q$ and $\qform{\la} q_\alpha$ represent a common element of $k$.  Then $q - \qform{\la} q_\alpha$ is an isotropic 16-dimensional form in $I^4 k$, hence it is hyperbolic.  Thus $\qform{\la} q_\alpha$ is represented by the trivial class in $H^1(k, SO(q))$; it follows that $\alpha$ is in the image of the map $H^1(k, Z(G_1)) \ra H^1(k, G_1)$ and that $G'_1$ is isomorphic to $G_1$.  Applying Tits's Witt-type theorem again, we find that $G$ and $G'$ are isomorphic.
\end{proof}

\begin{rmk}[$\chr k \ne 2$]
Let $G$ be as in Prop.~\ref{class.prop} and let $\gamma$ be the 3-Pfister form corresponding to $a(G)$.  We claim that the group $G_1$ in the proof of the proposition is isomorphic to $\Spin(\gamma^K)$, where $\gamma^K$ denotes the $K$-associate of $\gamma$ as defined in \cite[p.~499]{KMRT}.  In particular, \emph{the semisimple anisotropic kernel of $G$ is the semisimple anisotropic kernel of $\Spin(\gamma^K)$.}  To prove the claim, let $G'$ be the quasi-split strongly inner form of $G$, so $G'_1$ is the spin group of a quasi-split quadratic form $q$.  Fix $\alpha \in H^1(k, \Spin(q))$ such that $q_\alpha$ is isomorphic to $\gamma^K$.  The twisted group $G'_\alpha$ is as in the first three lines of the table and $a(G'_\alpha)$ equals the Arason invariant of $q_\alpha - q$.  This form is Witt-equivalent to $\qform{\delta} \gamma$ for $\delta \in \kx$ such that $K = k(\sqrt{\delta})$, hence $a(G'_\alpha)$ equals $a(G)$ and $G'_\alpha$ is isomorphic to $G$ by Prop.~\ref{class.prop}.  Since the semisimple anisotropic kernels of $G'_\alpha$ and $\Spin(q_\alpha)$ agree, the claim follows.
\end{rmk}

The following remarks on the proposition make some forward references, but we will not refer to them elsewhere in the paper.  

\begin{rmk}
The hypothesis that the indexes of \emph{both} $G$ and $G'$ are in the first three rows of the table is crucial.  For example, take $G$ to be the real Lie group EIII and let $G'$ be the compact real Lie group of type $E_6$.  The index of $G$ is in the second row of the table, but  $G'$ is anisotropic.  Nonetheless, $a(G)$ and $a(G')$ both equal the unique nonzero element of $H^3(\R, 2)$, as can be seen by combining \ref{rGAK} and \cite[p.~120]{Jac:ex}.
\end{rmk}

\begin{rmk}
Given a symbol $\gamma \in H^3(k, 2)$, there is a unique corresponding octonion algebra $C$.  The index of the group $G := G(\hrm_3(C, \qform{1, 1, -1}), K)$ appears in the first three rows of Table \ref{iso.table} by the main theorem for every $K$, and $a(G) = \gamma$ by \ref{rGAK}.  Combined with Propositions \ref{trivial.prop} and \ref{class.prop}, we conclude that \emph{every group whose index appears in the first three rows of Table \ref{iso.table} is isomorphic to $G(\hrm_3(C, \qform{1, 1, -1}), K)$ for some octonion algebra $C$ and some $K$.}  This is approximately the content of \cite[3.3]{Veld:unit1} and \cite[3.2]{Veld:unit2}.

Assuming the main theorem, we can rephrase the conclusion above as: \emph{Conditions (1)--(3) in Th.~\ref{MT} are equivalent to}
\begin{equation}
\parbox{4in}{\emph{$G(A, K)$ is isomorphic to $G(A', K)$ for some Albert $k$-algebra $A'$ with nonzero nilpotents.}}
\end{equation}
\end{rmk}

\section{Albert algebra reminders}  \label{remind}

\begin{borel}{Arbitrary Albert algebras} \label{ARALAL} 
Albert
algebras are Jordan algebras of degree 3 and hence may all be
obtained from cubic forms with adjoint and base point in the sense
of \cite{McC:FST}. More specifically, given an Albert algebra
$A$ over $k$, there exist a cubic form $N:A \rightarrow k$ (the
\emph{norm}) and a quadratic map $\sharp:A \rightarrow A, x
\mapsto x^\sharp$ (the \emph{adjoint}) which, together with the
unit element $1 \in A$, satisfy the relations $N(1) = 1, 1^\sharp
= 1$,
\begin{align}
\label{adj} x^{\sharp\sharp} &= N(x)x\,, \\
(DN)(x)y &= T(x^\sharp,y)\,, \notag \\
\label{unt} 1 \times x &= T(x)1 - x
\end{align}
in all scalar extensions, where $T = -(D^2\,\mbox{log}\,N)(1):A
\times A \rightarrow k$ (the \emph{trace form}) stands for the
logarithmic hessian of $N$ at $1$, $T(x) = T(x,1)$ and $x \times y
= (x + y)^\sharp - x^\sharp - y^\sharp$ is the bilinearization of
the adjoint. The $U$-operator of $A$ is then given by the formula
\begin{align}
\label{uop} U_xy = T(x,y)x - x^\sharp \times y\,.
\end{align}
The \emph{quadratic trace}, defined by $S(x) := T(x^\sharp)$,
is a quadratic form with bilinearization
\begin{align}
\label{bquatr} S(x,y) = T(x)T(y) - T(x,y)\,.
\end{align}
It relates to the adjoint by the formula
\begin{align}
\label{ads} x^\sharp &= x^2 - T(x)x + S(x)1\,,
\end{align}
where, as in \cite[1.5]{Jac:J}, the powers of $x \in A$ are
defined by $x^0 = 1$, $x^1 = 1$, and $x^{n+2} = U_xx^n$ for $n \ge 0$. For future use, we recall the following identities from
\cite{McC:FST}.
\begin{align}
\label{pad} x^\sharp \times (x \times y) &= N(x)y +
T(x^\sharp,y)x\,, \\
\label{eul} T(x^\sharp,x) &= 3N(x)\,, \\
\label{pau} x^\sharp \times x &= [S(x)T(x) - N(x)]1 - S(x)x -
T(x)x^\sharp\,, 
\end{align}
and
\begin{multline} \label{ppad} 
x^\sharp \times (y \times z) + (x \times y) \times (x
\times z) = T(x^\sharp,y)z + T(x^\sharp,z)y + T(y \times z,x)x\,.
\end{multline}
 Just as in general Jordan rings, the \emph{Jordan triple
product} derives from the $U$-operator \eqref{uop} through
linearization:
\begin{align} \label{jtp} 
\{x,y,z\} &:= U_{x,z}y := \big(U_{x+y} - U_x -
U_y\big)z \\
&= T(x,y)z + T(y,z)x - (z \times x) \times y\,. \notag
\end{align}
Finally, the generic minimum polynomial of $x$ in the sense of
\cite{JacKatz} is
\begin{align}
\label{genmin} m_x(\textbf{t}) = \textbf{t}^3 - T(x)\textbf{t}^2 +
S(x)\textbf{t} - N(x) \in k[\textbf{t}]\,,
\end{align}
so by \cite[p. 219]{JacKatz} we have
\begin{align}
\label{mineq} x^3 - T(x)x^2 + S(x)x - N(x)1 = 0 = x^4 - T(x)x^3 +
S(x)x^2 - N(x)x\,,
\end{align}
where the second equation is a \emph{trivial} consequence of the
first only for char $k \neq 2$ because then we are dealing with
\emph{linear} Jordan algebras. An element $x \in A$ is invertible
if and only if $N(x) \neq 0$. At the other extreme, $x \in A$ is
said to be \emph{singular} if $x^\sharp = 0 \neq x$. As an ad-hoc
definition, singular idempotents will be called \emph{primitive}.

The following lemma collects well-known properties of the Peirce
decomposition relative to primitive idempotents, using the
labelling of \cite{Loos}. We refer to \cite[Lemma
1.5]{Faulk:oct}, \cite[(28), (31), Lemma 2]{Racine:note}, and
\cite[Lemma 5.3d]{PR:mod2} for details. 

\begin{lem} \label{L-PEIRCEDEC}
Let $e \in
A$ be a primitive idempotent and put $f = 1- e$.
\begin{enumerate}
\item The Peirce components of $A$ relative to $e$ are
described by the relations
\begin{align*}
A_2(e) &= ke\,, \\
A_1(e) &= \{x \in A \mid T(x) = 0, e \times x = 0\}\,, \\
A_0(e) &= \{x \in A \mid e \times x = T(x)f - x\}\,.
\end{align*}
\item The restriction $S_e$ of $S$ to $A_0(e)$ is a quadratic
form with base point $f$ over $k$ whose associated Jordan algebra
agrees with $A_0(e)$ as a subalgebra of
$A$.
\item $x^\sharp = S(x)e$ for all $x \in A_0(e)$. \hfill $\qed$
\end{enumerate}
\end{lem}
\end{borel}

\begin{borel}{Reduced Albert algebras} \label{REDALAL} 
An Albert $k$-algebra $A$ is \emph{reduced} if it is not division, i.e., if it contains nonzero elements that are not invertible.  Every such algebra is isomorphic to a Jordan algebra $\hrm_3(C, \G)$ as defined in the following paragraph.

Let $C$ be an octonion (or Cayley) $k$-algebra (see \cite[\S33.C]{KMRT} for the
definition and elementary properties) and fix a diagonal matrix $\G := \qform{\gamma_1, \gamma_2, \gamma_3} \in \GL_3(k)$.  We write 
$\hrm_3(C, \G)$ for the vector space of 
$3$-by-$3$ matrices $x$ with entries in $C$ that are
$\Gamma$-hermitian ($x = \Gamma^{-1}\overline{x}^t\Gamma$) and
have scalars down the diagonal.  It is spanned by the diagonal unit vectors
$e_{ii}$ ($1 \leq i \leq 3$) and by the hermitian matrix units
$x_i[jl] := \gamma_lx_ie_{jl} +
\gamma_j\overline{x_i}e_{lj}\,\,(x_i \in C)$, where $(ijl)$ here
and in the sequel always varies over the cyclic permutations of
(123). It has a Jordan algebra structure derived from a cubic
form with adjoint and base point as in \ref{ARALAL}: Writing $N_C$
for the norm, $T_C$ for the trace and $x \mapsto \overline{x}$ for
the conjugation of $C$, we let
\begin{align*}
x = \sum \alpha_i e_{ii} + \sum x_i[jl]\,,\,y = \sum \beta_i
e_{ii} + \sum y_i[jl]\,, &&\text{($\alpha_i,\beta_i \in k, x_i,y_i
\in C$)}
\end{align*}
be arbitrary elements of $A$ and set
\begin{align}
\label{NORM} N(x) &= \alpha_1\alpha_2\alpha_3 - \sum
\gamma_j\gamma_l\alpha_iN_C(x_i) +
\gamma_1\gamma_2\gamma_3T_C(x_1x_2x_3)\,, \\
\label{ADJ} x^\sharp &= \sum \big(\alpha_j\alpha_l -
\gamma_j\gamma_lN_C(x_i)\big)e_{ii} + \sum
(\gamma_i\overline{x_jx_l} - \alpha_ix_i)[jl]\,,
\end{align}
as well as $1 = \sum e_{ii}$. Then $(N,\sharp,1)$ is a cubic form
with adjoint and base point whose associated trace form is
\begin{align}
\label{BILTR} T(x,y) = \sum \alpha_i\beta_i + \sum
\gamma_j\gamma_lN_C(x_i,y_i)\,.
\end{align}
The Jordan algebra structure on $\hrm_3(C, \G)$ is defined to be the one associated with $(N, \sharp, 1)$.

Specializing $y$ to $1$ in \eqref{BILTR} yields
\begin{align}
\label{LINTR} T(x) = \sum \alpha_i \,,
\end{align}
whereas linearizing \eqref{ADJ} leads to the relation
\begin{align} \label{POLADJ}
x \times y = & \sum \big(\alpha_j\beta_l +
\beta_j\alpha_l - \gamma_j\gamma_lN_C(x_i,y_i)\big)e_{ii} \\
& + \sum
(\gamma_i\overline{x_jy_l + y_jx_l} - \alpha_iy_i -
\beta_ix_i)[jl]\,. \notag
\end{align}
Furthermore, the quadratic trace $S$ of $A$ and its polarization by \eqref{ADJ},
\eqref{LINTR}, \eqref{POLADJ} have the form
\begin{align}
\label{QUADTRMAT} S(x) &= \sum \big(\alpha_j\alpha_l -
\gamma_j\gamma_lN_C(x_i)\big)\,, \\
\notag S(x,y) &= \sum \big(\alpha_j\beta_l +
\beta_j\alpha_l - \gamma_j\gamma_lN_C(x_i,y_i)\big)\,.
\end{align}
The triple $F := (e_{11},e_{22},e_{33})$ is a \emph{frame}, i.e., a complete
orthogonal system of primitive idempotents, in $A$ with Peirce
components $J_{ii}(F) = Re_{ii}\,,\,J_{jl}(F) = C[jl]$. Comparing
this with \eqref{QUADTRMAT}, we obtain a natural isometry
\begin{align}
\label{QUADTRPEDECOMP} S\vert_{J_{jl}(F)} \cong \langle
-\gamma_j\gamma_l \rangle\,.\,N_C\,.
\end{align}
\end{borel}

\begin{borel}{The split Albert algebra} 
Of
particular interest is the \emph{split} Albert algebra $A^d =
\hrm_3(C^d,1)$, where $C^d$ stands for
the split octonion algebra of Zorn vector matrices over $k$
\cite[Ch.~VIII, Exercise 5]{KMRT} and 1 is the 3-by-3
unit matrix. Albert algebras are exactly the $k$-forms of $A^d$,
i.e., they become isomorphic to $A^d$ over the separable closure
of $k$. For any Albert algebra $A$ over $k$ as in \ref{ARALAL},
the equation $N = 0$ defines a singular hypersurface in $\mathbb{P}(A)$,
its singular locus being given by the system of quadratic
equations $x^\sharp = 0$ in $\mathbb{P}(A)$. 
\end{borel}

\begin{borel}{The invariants $f_3$ and $f_5$} \label{INVMOD2}
Letting $A$ be a reduced Albert algebra over $k$ as in
\ref{REDALAL}, we briefly describe the cohomological mod-2
invariants of $A$ in a characteristic-free manner; for char $k
\neq 2$, see \cite[pp.~49--52]{GMS}. Since $C$ is
uniquely determined by $A$ up to isomorphism \cite[Th.~1.8]{Faulk:oct}, 
so is its norm $N_C$, which, being a $3$-fold
Pfister form, gives rise to the first mod-2 invariant $f_3(A) \in
H^3(k,2)$ of $A$. On the other hand, $\Gamma$ is not uniquely
determined by $A$ since, for example, it may be multiplied by
nonzero scalars and its components may be permuted arbitrarily as
well as multiplied by nonzero square factors without changing the
isomorphism class of $A$. This allows us to assume $\gamma_2 = 1$.
Then, as in \cite[2.1, 4.1]{Ptr:struct}, we may consider the 5-fold
Pfister form $\pform{-\gamma_1,-\gamma_3}
\, N_C$, giving rise to the second mod-2 invariant
\begin{align}
\label{MOD25} f_5(A) = f_3(A)\cdot (-\gamma_1) \cdot (-\gamma_3) \in
H^5(k,2)
\end{align}
of $A$. Translating Racine's characteristic-free version
\cite[Th.~3]{Racine:note} of Springer's classical result
\cite[p.~381, Th.~6]{Jac:J} to the cohomological setting, it
follows in all characteristics that reduced Albert algebras are
classified by their mod-2 invariants $f_3,f_5$.
\end{borel}

\begin{borel}{Nilpotent elements} \label{NILEL}  
Let $A$ be an
Albert algebra over $k$. An element $x \in A$ is said to be
\emph{nilpotent} if $x^n = 0$ for some $n \in \N$. Combining
\cite[p.~222, Th.~2(vi)]{JacKatz} with \eqref{genmin}, we conclude
that
\begin{align}
\label{CHARNIL} &&\text{$x \in A$ is nilpotent if and only if
$T(x) = S(x) = N(x) = 0$.}
\end{align}
In this case $x^3 = x^4 = 0$ by \eqref{mineq}. Hence $A$ contains
nonzero nilpotents if and only if  $x^2 = 0$ for some nonzero
element $x \in A$. We also conclude from \eqref{ads},
\eqref{CHARNIL} that $x \in A$ satisfies $x^2 = 0$ if and only if
$x^\sharp = 0$ and $T(x) = 0$. Finally, by \cite[4.4]{Ptr:struct},
an Albert algebra contains nonzero nilpotent elements if and only
if it is reduced and $f_5(A) = 0$. 
\end{borel}


In order to bring Jordan-theoretic techniques to bear in the proof of Theorem \ref{MT}, we will translate hypothesis (1) into a condition on the Albert algebra $A$.  To do so, we extend the results of \cite[\S13]{CG} to the case where $k$ is arbitrary.

\begin{defn}
A subspace $X$ in an Albert algebra $A$ is an \emph{inner ideal} if $U_x A \subseteq X$ for all $x \in X$.  A subspace $X$ of $A$ is \emph{singular} if $x^\sharp = 0$ for all $x \in X$.  A subspace $X$ of $A$ is a \emph{hyperline} if it is of the form $x \times A$ for some nonzero $x$ with $x^\sharp = 0$.
\end{defn}

The nonzero, proper inner ideals of an Albert algebra $A$ are the singular subspaces (of dimensions 1 through 6) and the hyperlines (of dimension 10), see \cite[p.~467]{McC:inn} and \cite[Th.~2]{Racine:point}.

For an inner ideal $X$ of $A$, we define $\psi(X)$ to be the set of all $y \in A$ satisfying
\begin{equation} \label{psi1}
\{ X, y, A \} \subseteq X
\end{equation}
and
\begin{equation} \label{psi2}
U_A U_y X \subseteq X.
\end{equation}

The definition of $\psi$ in \cite{CG} consisted of only condition \eqref{psi1}.
Condition \eqref{psi2} was suggested by Erhard Neher, for the purposes of including characteristic 2.  He points out that when $k$ has characteristic not 2, \eqref{psi2} follows from \eqref{psi1} by applying the identity JP13.  (The notation JP$xx$ refers to the Jordan pair identities as numbered in \cite{Loos}.)

\begin{lem}[Neher]
If $X$ is an inner ideal of $A$, then $\psi(X)$ is also an inner ideal.
\end{lem}

\begin{proof}
Clearly $\psi(X)$ is closed under scalar multiplication.  Let $y, z$ be elements of $\psi(X)$.  To prove that $y + z$ is in $\psi(X)$, it suffices to show that $U_A U_{y,z} X \subseteq X$, which follows from JP13:
\[
U_a U_{y,z} x = \{ a, y, \{ a, z, x \} \} - \{ U_a y, z, x \} \in X.
\]

Fix $y \in \psi(X)$ and $a \in A$.  By JP7, we have
\[
\{ X, U_y a, A \} \subseteq \{ \{ a, y, X\}, y, A \} + \{ a, U_y X, A \} \subseteq X,
\]
So $U_y a$ satisfies condition \eqref{psi1}.  Also, by JP3, 
\[
U_A U_{U_y a} X = U_A U_y U_a U_y X \subseteq X,
\]
so $U_y a$ satisfies condition \eqref{psi2}.
\end{proof}

We now check that the descriptions of $\psi(X)$ for various $X$ given in \cite[\S13]{CG} are valid in all characteristics.  Consider first the case where $X$ is 1-dimensional singular, i.e., an ``$\alpha_1$-space'' in the language of \cite{CG}.  For $y \in X \times A$, we have $\{ X, y, A \} \subseteq X$ by \cite[13.8]{CG}.  The argument in the last paragraph of the proof of \cite[Prop.~13.19]{CG} shows that $U_y X$ is zero.  Therefore, $\psi(X)$ contains $X \times A$, and as in \cite[13.6]{CG}---using that $\psi(X)$ is an inner ideal---we conclude that $\psi(X)$ equals $X \times A$.

Similar arguments easily extend the proofs of 13.9--13.19 of \cite{CG} to include the case where $k$ has arbitrary characteristic.  (See also Remarks 5.5 and 13.22 of that paper.)  For every $x$ in a proper inner ideal $X$ and every $y \in \psi(X)$, we find \emph{a posteriori} that $U_y x$ is zero.  That is, \eqref{psi2} is superfluous.

\begin{borel*} \label{CG.results}
In summary, for a proper, nonzero inner ideal $X$ in $A$, the set $\psi(X)$ of $y \in A$ satisfying \eqref{psi1} is an inner ideal of $A$ whose dimension is given by the following table.  The 5-dimensional inner ideals come in two flavors; we write ``$5'$" to indicate a 5-dimensional maximal singular subspace.
\begin{center}
\begin{tabular}{c|rrrrrr}
$\dim X$ & 1 & 2 & 3 & $5'$&6&10 \\ \hline
$\dim \psi(X)$ & 10 & $5'$& 3 & 2&6 & 1
\end{tabular}
\end{center}
For $X$ of dimension $\ne 6$, the space $\psi(X)$ is the set of all $y \in A$ such that $\{ X, y, A \}$ is zero.  For $X$ of dimension 6, $\psi(X)$ is a 6-dimensional inner ideal such that $\{ X, \psi(X), A \}$ equals $X$.
\end{borel*}

\section{The groups $G(A, K)$} \label{hjp}

This section defines the groups $G(A, K)$ using 
\emph{hermitian Jordan triples}.  These triples were previously studied in an analytic
setting \cite[2.9]{Loos:irv}; translating to a purely algebraic situation in a natural way, we arrive at the following concept.

\begin{defn} \label{triple}
A \emph{hermitian Jordan triple} over $k$ is a triple $(K, V, P)$ consisting of a quadratic \'etale $k$-algebra $K$ (with conjugation $\iota$), a free $K$-module $V$ of finite rank, and a quadratic map 
$P \!: V \ra \Hom_K(\iV, V)$, where $\iV$ is the $K$-module $V$ with scalar multiplication twisted by $\iota$, such that $(V, P)$ is an ordinary Jordan triple over $k$ in the sense of, e.g., \cite[1.13]{Loos}.  In particular, $P_v \!: V \ra V$ is an $\iota$-semilinear map depending $K$-quadratically on $v \in V$.

A homomorphism $(K, V, P) \ra (K', V', P')$ of hermitian Jordan triples is a pair $(\phi, h)$ such that $\phi \!: K \iso K'$ is a $k$-isomorphism and $h \!: V \ra V'$ is a $\phi$-semilinear map 
such that $h(P_x y) = P_{h(x)}h(y)$ for all $x, y \in V$.
\end{defn}

It is easy to see that Jordan pairs are basically the same as hermitian Jordan triples $(K, V, P)$ with $K$ split.

\begin{eg}
Starting with a Jordan $k$-algebra $J$ and a quadratic \'etale $k$-algebra $K$ with conjugation $\iota$, we obtain a hermitian Jordan triple as output.  Namely, take $V = J \ot_k K$ and define $P$ by
\[
P_x y := U_x \, (\Id_J \ot \iota) y.
\]
We denote this hermitian Jordan triple by $\T(J, K)$.

A \emph{hermitian Albert triple} is a triple $\T(A, K)$, where $A$ is an Albert $k$-algebra.
We call the triple $\Td := \T(\Ad, k \times k)$ the \emph{split Albert triple} or simply the \emph{split triple}.
\end{eg}

We will now explicitly describe the automorphism group of a hermitian Jordan triple $\T(J, K)$.
Following \cite[1.7, p.~1.23]{Jac:Ark},
the \emph{structure group} of $J$, denoted by $\Str(J)$, is the
subgroup of $\GL(J)$ consisting of all bijective linear maps
$g \!: J \rightarrow J$ that may be viewed as isomorphisms from $J$
onto an appropriate isotope or, equivalently, that have the
property that there exists a bijective linear map $g^\dag \!:J
\rightarrow J$ satisfying the relation $U_{g(x)} =
gU_x(g^\dag)^{-1}$ for all $x \in J$; obviously it is a group scheme. The assignment $g \mapsto
g^\dag$ is an order 2 automorphism of the structure
group of $J$.

\begin{lem}
The group of $k$-automorphisms of $\T(J, K)$ is generated by the element $(\iota, \Id_J \ot \iota)$ of order 2 and by the elements $(\Id_K, g)$ for $g \in \Str(J)(K)$ satisfying $g^\dag = (\Id_J \ot \iota) g (\Id_J \ot \iota)$.
\end{lem}

\begin{proof}
$(\iota, \Id_J \ot \iota)$ is clearly an automorphism of $\T(J, K)$.  Conversely, multiplying any automorphism $(\phi, h)$ of $\T(J, K)$ by $(\iota, \Id_J \ot \iota)$ if necessary, we are allowed to assume that $\phi$ is the identity on $K$.  The equation $h(P_x y) = P_{h(x)} h(y)$ is equivalent to $U_{h(x)} = h U_x (\Id_J \ot \iota) h^{-1} (\Id_J \ot \iota)$, which gives the claim.
\end{proof}

\begin{borel*} \label{aut.grp}
Since the trace form of an Albert algebra $A$ is a
non-degenerate symmetric bilinear form, it follows from \cite[p.~502]{McC:FST} that the structure group of $A$ agrees with its group of norm
similarities; viewed as a group scheme, it is a reductive
algebraic group with center of rank 1 and its semisimple part is simply connected of type $E_6$ \cite[14.2]{Sp:jord}.

Below, we only consider hermitian Jordan triples $\T = (V, K, P)$ that become isomorphic to $\Td$ over a separable closure $\ksep$ of $k$.  We write $G(\T)$ for the semisimple part of the identity component of $\aut(\T)$; when $\T$ is of the form $\T(A, K)$, we will write $G(A, K)$ for short.  This group is simply connected of type $E_6$ because it is isomorphic to the semisimple part of $\Str(\Ad)$ over $\ksep$.   It acts on the vector space underlying $\T$, and over a separable closure this representation is isomorphic to a direct sum of the usual representation of $E_6$ on $\Jd$ and its contragradient.  Since the direct sum of these two representations is defined over $k$, the group has trivial Tits algebras.  Also, $G(\T)$ is of type $\dE$ if and only if $K$ is a field.
\end{borel*}

\begin{borel}{Relation with Tits's construction} 
We now observe that the Lie algebra $\g$ of $G(A, K)$ is precisely the one obtained from $A$ and $K$ using Tits's construction from \cite{Ti:const}.  For the duration of this subsection, we assume that $k$ has characteristic $\ne 2, 3$.  (Tits makes this assumption in his paper, so it is harmless for our purposes.)  Suppose first that $K$ is $k \times k$.  From the first paragraph of \ref{aut.grp}, we see that $\g$ consists of those linear transformations on $A$ that leave the cubic form $N$ ``Lie invariant''.  
By
\cite[p.~189, Th.~5]{Jac:J1}, this is the same as the Lie algebra $R_0(A) \oplus D(A)$, where $R_0(A)$ is the vector space generated by the transformations ``right multiplication by a trace zero element of $A$'' and $D(A)$ denotes the derivation algebra of $A$; this is what one obtains from Tits's construction.

Now consider the case where $K$ is a field.  We may view the vector space $A \ot K$ underlying $\T(A, K)$ as the subspace of $A_K \times A_K$ fixed by the map $(x, y) \mapsto (\iota y, \iota x)$.  In this way, we can view $G(A, K)$ as the group whose $k$-points are those $(f, f^\dag) \in G(A_K, K \times K)(K)$ such that 
$\iota f^\dag \iota = f$.  The differential of $f \mapsto f^\dag$ is $d \mapsto -d^*$.  It follows that the Lie algebra $\g$ is the subalgebra of $R_0(A_K) \oplus D(A_K)$ consisting of elements fixed by the map $d \mapsto -\iota d^* \iota$.  When $d$ is right multiplication by an element of $A_K$, we have $d^* = d$ because the bilinear form $T$ is associative \cite[p.~227, Cor.~4]{Jac:J}.  The derivation algebra $D(A_K)$ is spanned by commutators of right multiplications, hence is fixed elementwise by the map $d \mapsto -d^*$.  Consequently, the Lie algebra of $G(A, K)$ is identified with $\sqrt{\la} \,R_0(A) \oplus D(A)$, where $\la \in \kx$ satisfies $K = k(\sqrt{\la})$.  This is the Lie algebra constructed by Tits.  (Jacobson \cite{Jac:ex} and Ferrar \cite{Ferr:E6} denote this Lie algebra by $\mathscr{E}_6(A)_\la$ and $\mathscr{L}(A)_\la$ respectively.)
\end{borel}

By Galois descent, the map $\T \mapsto G(\T)$ produces (up to isogeny) all groups of type $E_6$ with trivial Tits algebras over every field $k$.  But over some fields, it suffices to consider only the triples of the form $\T(A, K)$.

\begin{eg}
Every group of type $\oE$ with trivial Tits algebras can be realized as the group of linear transformations on an Albert algebra $A$ preserving the cubic norm, hence is of the form $G(A, k \times k)$.  (In characteristic $\ne 2, 3$, this follows from \cite[3.4]{G:rinv}.)  We now list a few examples of fields $k$ such that every group of type $\dE$ with trivial Tits algebras is---up to isogeny---of the form $G(A, K)$.
\begin{enumerate}
\item When $k$ is a local field \cite{Kn:p2} or a finite field, every group of type $\dE$ is quasi-split.

\item For $k$ the real numbers, up to isogeny the three groups of type $\dE$ are of the form $G(A, \C)$ as $A$ varies over the three isomorphism classes of Albert $\R$-algebras \cite[pp.~119, 120]{Jac:ex}.

\item When $k$ is a number field, as can be seen by using the Hasse Principle, cf.~\cite[3.2, 6.4]{Ferr:E62}.
\end{enumerate}
\end{eg}

For $f$ an automorphism of an Albert algebra $A$, we have $f^\dag = f$, so $\aut(A)$ is a subgroup of $G(A, K)$ for all $K$.

\begin{eg}
The group $G(A^d, K)$ contains the split group $\aut(A^d)$ of type $F_4$, hence has $k$-rank at least 4.  The classification of possible indexes gives that $G(A^d, K)$ is quasi-split for all $K$, including the case $K = k \times k$.
\end{eg}

\begin{borel}{Rost invariant of $G(A, K)$} \label{rGAK}
The inclusion of $\aut(A)$ in $G(A, K)$ has Rost multiplier one \cite[p.~194]{Dynk:ssub} in the language of \cite{GMS}, hence the composition
    \[
    \begin{CD}
    H^1(k, \aut(A)) @>>> H^1(k, G(A, K)) @>{r_{G(A, K)}}>> H^3(k, n_{G(A, K)})
    \end{CD}
    \]
agrees with the Rost invariant relative to $\aut(A)$.  In the definition of $a(G)$ at the start of \S\ref{iso}, one can take the image of $\Ad$ as $\eta$ by the preceding example.  We find:
   \[
   a(G(A, K)) = -r_{F_4}(A) \quad \in H^3(k, 6),
      \]
where the right side denotes the negative of the usual Rost invariant of the Albert  algebra $A$.  In particular, the mod-2 portion of $a(G(A, K))$ is the symbol $f_3(A)$.
\end{borel}

If one knows that a group $G(A, K)$ is isotropic, then combining the previous paragraph and Prop.~\ref{trivial.prop} gives the index of $G(A, K)$.

\begin{eg} \label{bad.eg}
Fix a group $G$ of type $\dE$ with trivial Tits algebras whose index is from the bottom row of Table \ref{iso.table}.  By Prop.~\ref{trivial.prop}, $a(G)$ is not a symbol, hence $G$ is not of the form $G(A, K)$ for any Albert algebra $A$. 

We remark that the use of the Rost invariant simplifies this example dramatically, compare \cite[pp.~64, 65]{Ferr:E6}.
\end{eg}

\section{Homogeneous projective varieties} \label{homo}

In this section, we relate the Tits index of a group $G(\T)$ to inner ideals in the triple $\T$.  This is done by describing the $k$-points on the homogeneous projective varieties for the groups.

A projective variety $Z$ is \emph{homogeneous} for $G(\T)$ if $G(\T)$ acts on $Z$ and $G(\T)(\ksep)$ acts transitively on $Z(\ksep)$, for $\ksep$ a separable closure of $k$.  There is a bijection between the collection of such varieties defined over $k$ (up to an obvious notion of equivalence) and Galois-stable subsets of vertices in the Dynkin diagram \cite[6.4]{BoTi}.  There are two common ways of normalizing the bijection: we normalize so that the trivial variety corresponds to the empty set and the largest homogeneous variety (i.e., the variety of Borel subgroups) corresponds to the full Dynkin diagram.

Write $\T = (K, V, P)$.  
A $K$-submodule $X$ of $V$  is an \emph{inner ideal in $\T$} if $P_x V \subseteq X$ for every $x \in X$.  
We write $\br{\ }$ for the bilinearization of $P$, i.e., 
\[
\br{x, y, z} := P_{x + z} y - P_x y - P_z y.
\]
Since $P$ is a quadratic map and $P_x$ is $\iota$-semilinear for all $x \in V$, the ternary product $[\ ]$ is linear in the outer two slots and $\iota$-semilinear in the middle slot.

\begin{eg}\label{inn}
When $K$ is a field, the inner ideals of $\T(A, K)$ are the same as the inner ideals of the Albert $K$-algebra $A_K$.  To see this, note that for a $K$-submodule $X$ of $A_K$, we have $P_X A_K = U_X (\Id_A \ot \iota) A_K = U_X A_K$.

When $K$ is not a field, i.e., when $K$ is $k \times k$, the $K$-module underlying $\T(A, K)$ is $A \times A$.  A $K$-subspace $X$ of $A \times A$ is of the form $X_1 \times X_2$ for $X_i$ a $k$-subspace of $A$.  The same argument as in the previous paragraph gives that $X$ is an inner ideal of $\T(A, K)$ if and only if $X_1$ and $X_2$ are inner ideals of $A$.
\end{eg}

Below, we only consider inner ideals that are free $K$-modules, and by \emph{rank} we mean the rank as a $K$-module.  In the case where $K$ is a field, the previous example gives: The rank one (resp., 10) inner ideals of $\T(A,K)$ are of the form $Kx$ (resp., $x \times A_K$) where $x \in A_K$ is singular. 

\begin{prop} \label{flag}
The projective homogeneous variety for $G(\T)$ corresponding to the subset $S$ of the Dynkin diagram has $k$-points the inner ideals as in the table below.
\[
\begin{tabular}{cc}
$S$& $k$-points are inner ideals\\ \hline
$\{ \alpha_1, \alpha_6 \}$&$X \subset Y$ s.t.\ $\rank X = 1$, $\rank Y = 10$, and $\br{X, Y, V} = 0$ \\
$\{ \alpha_3, \alpha_5 \}$&$X \subset Y$ s.t.\ $\rank X = 2$, $\rank Y = 5$, and $\br{X, Y, V} = 0$\\
$\{ \alpha_4 \}$&$X$ s.t.\ $\rank X = 3$ and $\br{X, X, V} = 0$\\
$\{ \alpha_2 \}$&$X$ s.t.\ $\rank X = 6$ and $\br{X, X, V} = X$
\end{tabular}
\]
\end{prop}

With the information from the table, it should be no trouble to describe the homogeneous varieties corresponding to other opposition-stable subsets of the Dynkin diagram using the recipe in \cite[9.4]{CG}.  (The \emph{opposition involution} is the unique nonidentity automorphism of the Dynkin diagram.)

\begin{proof}
Since the claimed collections of subspaces form projective $k$-varieties on which $G(\T)$ acts, it remains only to show that the action is transitive on $\ksep$-points and that the stabilizer of a $\ksep$-point is a parabolic subgroup of the correct type.  Therefore, we assume that $k$ is separably closed---in particular that $\T$ is split, i.e., isomorphic to $\Td$---and note that $G(\T)$ is the split simply connected group of type $E_6$.  The projective homogeneous variety $Z_0$ associated with the split group and $S$ was described in sections 7 and 9 of \cite{CG}.  We prove the proposition by producing a $G(\T)(k)$-equivariant bijection $Z(k) \iso Z_0(k)$.

\smallskip

Let $X \subset Y$ be inner ideals as in the first line of the table.  As in Example \ref{inn}, we write $X = X_1 \times X_2$ and similarly for $Y$, where $X_i$ and $Y_i$ are inner ideals in $\Jd$.  The set $Z_0(k)$ consists of pairs $X' \subset Y'$ where $X'$ and $Y'$ are inner ideals in $\Jd$ of dimension 1 and 10 respectively.  Define $f \!: Z \ra Z_0$ by $f(X \subset Y) = (X_1 \subset Y_1)$; it is $G$-equivariant.

On the other hand, a pair $X' \subset Y'$ in $Z_0(k)$ is the image of $(X' , \psi(Y')) \subset (Y', \psi(X'))$ under $f$.  Indeed, $\psi(Y')$ and $\psi(X')$ are inner ideals of the  appropriate dimension by \ref{CG.results}, and $X' \subset Y'$ implies $\psi(Y') \subset \psi(X')$.

\smallskip

The second and third lines of the table follow by similar reasoning.

\smallskip

Now let $X$ be as in the last line of the table.  The $k$-points of $Z_0$ are the 6-dimensional singular subspaces of $\Jd$.  Again writing $X$ as $(X_1, X_2)$, we define $f$ by $f(X) = X_1$.
The equation $\br{X, X, V} = X$ is equivalent to the equations $\{ X_i, X_{i+1}, \Jd \} = X_i$ for $i = 1, 2$.  
Clearly, $X_{i+1}$ is a subset of $\psi(X_i)$.  But $\psi(X_i)$ also has dimension 6 by \ref{CG.results}, hence $X_{i+1} = \psi(X_i)$.  That is, $X_2$ is determined by $X_1$.  The conclusion now follows as in the case of the first line.
\end{proof}

We now use the description of the projective homogeneous varieties of $G(A, K)$ given in the proposition to give a concrete criterion for determining whether $G(A, K)$ is isotropic.

\begin{cor} \label{transl}
The group $G(A, K)$ is isotropic if and only if there is a nonzero $x \in A_K$ such that
   \[
   x^\sharp = 0 \quad \text{and} \quad x \in \iota(x) \times A_K.
   \]
\end{cor}

\begin{proof}
Combining \ref{rGAK} and Prop.~\ref{trivial.prop}, we find that $G(A, K)$ is isotropic if and only if the vertices $\alpha_1$ and $\alpha_6$ in its index are circled, i.e., if and only if the projective homogeneous variety corresponding to $\{ \alpha_1, \alpha_6 \}$ has a $k$-point \cite[6.4]{BoTi}.    Hence $G(A, K)$ is isotropic if and only if there are inner ideals $X \subset Y$ in $A_K$ of ranks 1 and 10 such that $[X, Y, A_K] = 0$.  

Suppose first that such inner ideals $X \subset Y$ exist.  We have $X = Kx$ for some nonzero $x \in A_K$ satisfying $x^\sharp = 0$ and
\[
\{ X, \iota(Y), A_K \} = [X, Y, A_K] = 0.
\]
Since \eqref{psi1} implies \eqref{psi2}, $\iota(Y)$ is contained in $\psi(X) = x \times A_K$.  Comparing ranks using \ref{CG.results}, we conclude that $X \subseteq Y = \iota(x) \times A_K$.

Conversely, suppose that $x \in A_K$ satisfies the conditions displayed in the corollary.  Then $X := Kx$ and $Y := \iota(x) \times A_K$ are inner ideals of $A_K$ of the desired ranks such that $X \subset Y$.  Also,
$\br{x, \iota(x) \times A_K, A_K} = \{ x, x \times A_K, A_K \}$, which is zero by \eqref{ppad} and \eqref{jtp}.
\end{proof}

\smallskip 

\begin{eg} \label{nil.eg2}
\emph{If $A$ has nonzero nilpotent elements, then $G(A, K)$ is isotropic.}  Indeed, $A$ has a nonzero element $x$ such that $x^\sharp = 0$ and $T(x) = 0$ by \ref{NILEL}.  Then $\iota(x) \times (-1) = x$ by \eqref{unt}.
\end{eg}

This example shows that in Theorem \ref{MT}, (4) implies (1).

We now resolve an open question from \cite{Veld:unit2}.

\begin{eg} \label{Veld}
Consider again a group $G$ as in Example \ref{bad.eg}, and write $K$ for the associated quadratic extension.  The group is of the form $G(\T)$ for a triple $\T = (V, K, P)$ where $\T$ contains an ``$\alpha_2$-space'', i.e., an inner ideal $X$ of $V$ such that $\dim X = 6$ and $\br{X, X, V} = X$.

Over $K$, $\T \ot K$ is isomorphic to a triple $\T(A, K \times K)$ for some Albert $K$-algebra $A$.  However, $A$ contains a 6-dimensional singular subspace, so it is split \cite[Th.~1]{Racine:point}. That is, $\T \ot K$ is isomorphic to the split triple $\T(A^d_K, K \times K)$.  (Alternatively, one can see this using Tits's list of possible indexes for groups of type $\oE$.)  By Galois descent, we may identify $V$ with 
 \[
 \{ (a, t \iota a) \mid a \in A^d_K \},
 \]
where $t$ is some $K$-linear transformation of $A^d_K$.  (Our $t\iota$ is Veldkamp's $T$.)  As in the proof of Prop.~\ref{flag}, we have $t \iota X = \psi(X)$, where we have identified $X$ with its first component in $A_K \times A_K$.  For $x \in X$, we have $T(x, t \iota x) = 0$ by \cite[13.19]{CG}. 

Veldkamp was concerned with geometries whose points were 1-dimensional singular $K$-subspaces of $A^d_K$.  In his language, the last sentence of the previous paragraph shows that $Kx$ is a weakly isotropic point when $x$ is nonzero.  Veldkamp asked on page 291 of \cite{Veld:unit2} if it is possible to have a weakly isotropic point and no strongly isotropic point.  We now observe that there is no strongly isotropic point in this example.  Suppose that $z \in A^d_K$ is nonzero and strongly isotropic, i.e., $\{ z, t \iota z, A^d_K \}$ is zero.  In that case, taking $Z = Kt \iota z$ and $Y = z \times A^d_K$ we find a pair of inner ideals of dimension 1 and 10 in $\T$ such that $Z \subset Y$ and 
$\br{Z, Y, V} = 0$.  But this contradicts the hypothesis on the index of $G(\T)$, so no such $z$ can exist, i.e., we have produced an explicit example of a situation where there is a weakly isotropic point and no strongly isotropic point.
\end{eg}

\section{Embeddings of $k \times K$} \label{emb}

In this section, we assume that $K$ is a separable quadratic field extension of $k$, and write $\iota$ for its nontrivial $k$-automorphism.  We fix a representative $\delta \in \kx$ for the discriminant of $K/k$.  (We take the naive definition of discriminant from, say, \cite{Lang}.  Consequently, in characteristic 2, we can and do take $\delta = 1$.)
The purpose of this (long) section is to prove the following proposition, which in turn shows that (2) implies (3) in the main theorem, see \ref{imp23}.
 
\begin{prop} \label{emb.prop}
Let  $A$ be an Albert $k$-algebra.  Then $k \times K$ is a subalgebra of $A$ if and only if $A$ is isomorphic to $\hrm_3(C, \qform{r, 1, \delta N_K(s) } )$ for some octonion algebra $C$ and elements $r \in \kx$ and $s \in \Kx$ such that $T_K(s) \ne 0$.
\end{prop}

The assumption that $K$ is a field is harmless.  If $K$ is $k \times k$, then we can take $\delta = 1$ and the proposition is easily seen to hold so long as $k$ is not the field with two elements.

\begin{borel}{\textsl{Proof that (2) implies (3)}} \label{imp23}  Before proceeding with the proof of the proposition, we first observe that it shows that (2) implies (3) in the main theorem.
Prop.~\ref{emb.prop} gives that $A$ is isomorphic to $\hrm_3(C, \qform{r, 1, \delta N_K(s)} )$.  That is, $f_3(A)$ is $[C]$ and $f_5(A)$ equals $[C] \cdot (- r) \cdot (-\delta N_K(s))$ by \eqref{MOD25}.  The 2-Pfister bilinear form $\gamma := \pform{-r, -\delta N_K(s)}$ satisfies (3).
\hfill$\qed$\end{borel}

\begin{borel}{\textsl{Proof that (4) implies (2)}} \label{imp42}  
The proposition also gives a proof that (4) implies (2) in the main theorem.    Specifically, if $A$ has a nilpotent, it is isomorphic to $\hrm_3(C, \qform{-1, 1, \delta N_K(s)})$ for every $\delta \in \kx$ and every $s \in \Kx$ with nonzero trace, since both algebras have the same invariants $f_3 = [C]$ and $f_5 = 0$.  Then $A$ contains $k \times K$ by the proposition.
\end{borel}

\begin{borel}{Proof of \ref{emb.prop}: the easy direction}  \label{easy}
We suppose that $A$ is $\hrm_3(C, \Gamma)$ for $\G$ as in the proposition and produce an explicit embedding of $k \times K$ in $A$.  Fix $t \in K^\times$ such that $t^2 = \delta$, which implies
\begin{align}
\label{DISC} T_K(t) = 0 \quad \text{and} \quad  N_K(t) = -\delta,
\end{align}
and define a map $\varphi \!: k \times K \rightarrow A$ by
\begin{align*}
\varphi(\alpha, a) := &
\alpha e_{11} + \\
&T_K(s)^{-1}\Big(N_K(s,a)e_{22} + N_K\big(s,\iota(a)\big)e_{33} +
t^{-1}\big(\iota(a) - a\big)1_C[23]\Big)
\end{align*}
for $\alpha \in k, a \in K$. Note that $t^{-1}(\iota(a) - a)$ is
in $k$ in all characteristics by \eqref{DISC}, i.e., the image of
$\varphi$ really is in $A$ and not just in $A_K$. Also, $\varphi$
sends 1 to 1. Hence it suffices to show that $\varphi$ preserves norms.  Plugging in to the explicit formula for the norm on $J$ from  \eqref{NORM}, we find 
\[
N \left[ \varphi(\alpha, a) \right] = \frac{\alpha}{T_K(s)^2} \left( 
N_K(s,a)N_K\big(s,\iota(a)\big) -
N_K(s)\big(\iota(a) - a\big)^2 \right).
\]
Writing out, for example, $N_K(s, a)$ as $s \iota(a) + \iota(s) a$, we find that
$N \left[ \varphi(\alpha, a) \right]$ is $\alpha N_K(a)$ as desired. $\hfill\qed$
\end{borel}

The converse implication of Prop.~\ref{emb.prop} takes a bit more work.
Fortunately, we can rely on the proof of \cite[Th.~3.2]{PR:Sp} to simplify the task. For the sake of clarity,
we still indicate the main steps of that proof insofar they are
relevant for our purposes. 

\begin{borel}{Proof of \ref{emb.prop}: the difficult direction}
Assume now
that $E := k \times K$ is a subalgebra of $A$ as in the statement of Prop.~\ref{emb.prop}.  We let $\iota$ act on $A_K
:= A \otimes K$ and $E_K := E \otimes K \subset A_K$ through the
second factor. Since $A$ inherits zero divisors from $E$, it is
reduced. We write $C$ for its coordinate algebra, forcing
$C_K := C \otimes K$ to be the coordinate algebra of
$A_K$. 

\smallskip 

\noindent\emph{Step 1. A coordinatization for $E_K$.} 
Since
$E_K$ is split but $E$ is not (because $K$ is a field by
hypothesis), there is a frame $(e_1,e_2,e_3)$ of $A_K$
satisfying
\begin{align*}
E_K = Ke_1 \oplus Ke_2 \oplus Ke_3\ \text{with}\ \iota(e_1) =
e_1\,,\,\iota(e_2) = e_3\,,\,\iota(e_3) = e_2\,.
\end{align*}
This setup will remain fixed till the end of the proof. It
implies
\[
E = \{\alpha e_1 + ae_2 + \iota(a)e_3\mid \alpha \in k, a \in
K\}\,. \vspace{4pt}
\]

\smallskip

\noindent\emph{Step 2. A coordinatization for $A_K$.} Following
\cite[Lemma 3.5]{PR:Sp}, we find a coordinatization of
$A_K$ having the form
\begin{align}
\label{COORD} A_K =
\hrm_3(C_K,
\qform{1,s,\iota(s)}) \quad \text{for some $s \in
K^\times$}
\end{align}
such that $e_i = e_{ii}$ for all $i$, the trace $T_K(s)$ is not zero, and 
there is an $\iota$-semilinear involution $\tau_K$ of
$C_K$ that commutes with the conjugation of $C _K$ and
satisfies the relation
\begin{align*}
\iota(x) = \iota(a_1)e_1 + \iota(a_3)e_2 + \iota(a_2)e_3 +
\tau_K(x_1)[23] + \tau_K(x_3)[31] + \tau_K(x_2)[12]
\end{align*}
for all $x = \sum a_ie_{ii} + \sum x_i[jl] \in A_K$. In particular,
\begin{equation} \label{DESC} 
A = \left\{ \left.
\parbox{2in}{
$\alpha e_1 + ae_2 + \iota(a)e_3$ \\
$+ x_1[23] 
+ x_2[31] + \tau_K(x_2)[12]$
}
\right| 
\,
\parbox{1.25in}{
$\alpha \in k, a \in K,$ \\
$x_1 \in C_K^{\tau_K}, x_2 \in C_K$}
\right\}\,,
\end{equation}
where $C_K^{\tau_K}$ denotes the subspace of $C_K$ consisting of elements fixed by $\tau_K$.

\smallskip

\noindent\emph{Step 3. A frame for $A$.} By \cite[Lemma
3.6]{PR:Sp}, $F = (d_1,d_2,d_3)$, where
\begin{align}
 d_1 &:= e_{11}\,,\notag \\
d_2 &:= T_K (s)^{-1}\big(se_{22} + \iota(s)e_{33} + 1[23]\big)
\label{FRA} \\
d_3 &:= T_K(s)^{-1}\big(\iota(s)e_{22} + se_{33} - 1[23]\big),
\notag
\end{align}
is a frame of $A$. Furthermore, the map
\begin{align}
\label{CO-DESC} \tau \!:C \ra C \ \text{defined by}\ x \mapsto \tau(x)
:= \overline{\tau_K(x)}\,,
\end{align}
is an $\iota$-semilinear automorphism of $C_K$, forcing $B := C^{\tau}_K$ to be a
composition algebra over $k$ such that $C_K = B \otimes K$.

\smallskip

\noindent\emph{Step 4. Peirce components of $A$.} As in \ref{easy}, we now fix
$t \in K^\times$ satisfying $t^2 = \delta$, hence \eqref{DISC}. Then the
following relations hold:
\begin{align}
\label{PEIRCE-12} A_{12}(F) &= \{x_2[31] + \overline{x_2}[12]\mid x_2 \in B\}\,, \\
\label{PEIRCE-31} A_{31}(F) &= \{(st)x_2[31] +
\iota(st)\overline{x_2}[12]\mid x_2 \in B\}\,.
\end{align}
While \eqref{PEIRCE-12} was established in \cite[Lemma 3.6b]{PR:Sp},
the proof of \eqref{PEIRCE-31} is a bit more
involved and runs as follows. Standard facts about Peirce
components and \eqref{DESC}, \eqref{FRA} imply
\begin{align*}
A_{31}(F) &= A_1(d_3) \cap A_1(d_1) = A_1(d_3) \cap A \cap
(C_K[31] + C_K[12]) \\
&= A_1(d_3) \cap \{x_2[31] + \tau_K(x_2)[12]\mid x_2 \in
C_K\}\,.
\end{align*}
For $x_2 \in C_K$, the element $x := x_2[31] +
\tau_K(x_2)[12]$ satisfies
\begin{align*}
T_K(s)(d_3 \times x) &= \big(\iota(s)e_{22} + se_{33} - 1[23]\big)
\times (x_2[31] + \tau_K(x_2)[12]) \\
&= -\big(s\overline{\tau_K(x_2)} + \iota(s)x_2\big)[31] -
\big(\iota(s)\overline{x_2} + s\tau_K(x_2)\big)[12]
\end{align*}
by \eqref{POLADJ} and \eqref{COORD},
and $x$ belongs to $A_1(d_3)$ if and only if this expression is
zero (Lemma \ref{L-PEIRCEDEC}a). Since $\iota(t) = -t$ by \eqref{DISC},
we may apply \eqref{CO-DESC} to conclude that this in turn is
equivalent to the equation $st\tau(x_2) = \iota(st)x_2$, i.e., to
$x_2 \in (st)B$.

\smallskip

\noindent\emph{Step 5. Conclusion of proof.} Using \eqref{QUADTRMAT} and
\eqref{COORD}, it is straightforward to check the relations
\begin{align}
\label{QUADTRPEIRCE-12} S(x_2[31] + \overline{x_2}[12]) &=
-T_K(s)N_B(x_2)\,, \\
\label{QUADTRPEIRCE-31} S\big((st)x_2[31] +
\iota(st)\overline{x_2}[12]\big) &= -T_K(s)N_K(s)\delta N_B(x_2)
\end{align}
for all $x_2 \in B$.  Following \cite[p.~1077]{McC:NAS}, we now choose a coordinatization of $A$ that has the form 
\begin{align*}
A = \hrm_3(C,\Gamma^0)\ \text{for $\Gamma^0 =
\qform{\gamma_1^0,1,\gamma_3^0} \in \GL_3(k)$,}
\end{align*}
and that matches $F$ with the diagonal frame of
$\hrm_3(C,\Gamma^0)$. Comparing  \eqref{QUADTRPEIRCE-12},
\eqref{QUADTRPEIRCE-31} with \eqref{QUADTRPEDECOMP}, we conclude
\begin{align}
\label{NORMPEIRCE-12} \qform{-\gamma_1^0} N_C &\cong
S\vert_{A_{12}(F)} \cong \qform{-T_K(s)} N_B\,, \\
\label{NORMPEIRCE-31} \qform{-\gamma_3^0\gamma_1^0} N_C
&\cong S\vert_{A_{31}(F)} \cong 
\qform{-T_K(s)N_K(s)\delta}N_B\,.
\end{align}
In particular, the composition algebras $B,C$ over $k$ have
similar norm forms and hence are isomorphic, allowing us to identify $B = C$ from now on. Then
\eqref{NORMPEIRCE-12} and \eqref{NORMPEIRCE-31}
imply that $A$ and $\hrm_3(C, \qform{T_K(s), 1, \delta N_K(s)})$ have the same mod-2 invariants, hence are isomorphic.
\hfill$\qed$
\end{borel}

\section{The case where $T(x) = 0$} \label{Tx0}

In this and the following section, we will prove that (1) implies (2) in the main theorem.  Suppose that (1) holds, so that by Cor.~\ref{transl} $A_K$ contains a nonzero element $x$ satisfying 
\begin{equation} \label{onestar}
x^\sharp = 0 \quad \text{and} \quad x \in \iota(x) \times A_K.
\end{equation}
In this section, we treat the case where $x$ satisfies \eqref{onestar} and has trace zero.  The next section treats the case where the trace of $x$ is nonzero.

\begin{prop} \label{Tx0.prop}
Let $A$ be an Albert $k$-algebra.
There is a nonzero element $x \in A_K$ satisfying
\[
x^\sharp = 0, \quad 
x \in \iota(x) \times A_K, \quad \text{and} \quad
T(x) = 0
\]
if and only if $A$ contains nonzero nilpotent elements.
\end{prop}

\begin{proof} 
We suppose that $A_K$ contains an $x$ as in the statement of the proposition; the other direction was treated in Example \ref{nil.eg2}.  If $K$ is equal to $k \times k$, then $A_K = A \times A$ and $x = (x_1, x_2)$ for $x_i \in A$ not both zero with $x_i^\sharp = 0$ and $T(x_i) = 0$.  Therefore one of the elements $x_1, x_2 \in A$ is a nonzero nilpotent by \ref{NILEL}.  We are left with the case where $K$ is a field.

For sake of contradiction, we further assume that $A$ has no nonzero nilpotents.  Write $K = k[d]$ where $d \in K$ has trace 1. Set
$\delta := N_K(d) \in k^\times$, so that $d^2 = d - \delta$.  Write
\[
x = y + dz \quad \text{with $y, z \in A$ not both zero.}
\]
Because $x^\sharp = 0$, we have:
\begin{equation} \label{yzz1}
y^\sharp = \delta z^\sharp \quad \text{and} \quad z^\sharp = -y \times z.
\end{equation}
Clearly, 
\begin{equation} \label{yzz2}
T(y) = T(z) = 0.
\end{equation}

We first argue that neither $y$ nor $z$ is invertible in $A$.  By \eqref{adj} we have:
\begin{equation} \label{rk1}
N(y)y = (y^\sharp)^\sharp = \delta^2 (z^\sharp)^\sharp = \delta^2N(z)z\,.
\end{equation}
On the other hand, comparing
\begin{align*}
y^\sharp \times (y \times z) &= N(y)z + T(y^\sharp,z)y &&\text{(by
\eqref{pad})} \\
&= N(y)z + \delta T(z^\sharp,z)y &&\text{(by \eqref{yzz1})} \\
&= N(y)z + 3\delta N(z)y &&\text{(by \eqref{eul})}
\end{align*}
with
\begin{align*}
y^\sharp \times(y \times z) &= -\delta z^\sharp \times z^\sharp
&&\text{(by \eqref{yzz1})} \\
&= -2\delta z^{\sharp\sharp} = -2\delta N(z)z\,, &&\text{(by
\eqref{adj})}
\end{align*}
we obtain
\begin{align}
\label{Nz3y} -2\delta N(z)z = N(y)z + 3\delta N(z)y \,.
\end{align}
Similarly,
\begin{align*}
\delta z^\sharp \times (z \times y) &= \delta N(z)y + \delta
T(z^\sharp,y)z &&\text{(by \eqref{pad})} \\
&= \delta N(z)y + T(y^\sharp,y)z &&\text{(by \eqref{yzz1})} \\
&= \delta N(z)y + 3N(y)z &&\text{(by \eqref{eul})}
\end{align*}
and
\begin{align*}
\delta z^\sharp \times (z \times y) &= -\delta z^\sharp \times
z^\sharp &&\text{(by \eqref{yzz1})} \\
&= -2\delta z^{\sharp\sharp} = -2\delta N(z)z &&\text{(by
\eqref{adj})}
\end{align*}
imply
\begin{align}
\label{Nzy3} -2\delta N(z)z = 3N(y)z + \delta N(z)y \,.
\end{align}
Subtracting \eqref{Nz3y} from \eqref{Nzy3}, we conclude $2N(y)z =
2\delta N(z)y$, which implies
\begin{align}
\label{Nyz} N(y)z = \delta N(z)y
\end{align}
for char $k \neq 2$, while this follows directly from \eqref{Nz3y}
for char $k = 2$. Thus \eqref{Nyz} holds in full generality. By
\eqref{rk1}, assuming that one of the elements $y,z$ is invertible,
both are, and
\begin{align*}
N(y)^2z^\sharp &= \big(N(y)z\big)^\sharp = \big(\delta
N(z)y\big)^\sharp &&\text{(by \eqref{Nyz})} \\
&= \delta^2N(z)^2y^\sharp = \delta^3N(z)^2z^\sharp &&\text{(by
\eqref{yzz1})}
\end{align*}
implies
\begin{align}
\label{NyN} N(y)^2 = \delta^3N(z)^2\,.
\end{align}
On the other hand,
\begin{align*}
\delta N(y)z^\sharp &= -\delta y \times [N(y)z] &&\text{(by
\eqref{yzz1})} \\
&= -\delta^2N(z)y \times y &&\text{(by \eqref{Nyz})} \\
&= -2\delta^2N(z)y^\sharp = -2\delta^3N(z)z^\sharp\,, &&\text{(by
\eqref{yzz1})}
\end{align*}
which yields
\begin{align}
\label{Ny2} N(y) = -2\delta^2N(z)\,.
\end{align}
Hence char $k \neq 2$, and comparing the square of \eqref{Ny2}
with \eqref{NyN} shows $\delta = \frac{1}{4}$. But then the
minimum polynomial of $d$ over $k$ becomes $X^2 - X + \frac{1}{4}
= (X - \frac{1}{2})^2$, contradicting the fact that $K$ is a
field. We have thus established the relation
\begin{align}
\label{N0} N(y) = N(z) = 0\,.
\end{align}

By assumption,
$A$ does not contain nonzero nilpotent elements.
Hence \eqref{yzz2}, \eqref{N0} imply that $S(y), S(z)$ cannot both
be zero since, otherwise, $y,z$ were both nilpotent by
\eqref{CHARNIL}. But since $S(y) = \delta S(z)$ by \eqref{yzz1}, we
conclude $S(y)$ and $S(z)$ are both nonzero. On the other hand, \eqref{N0}
combines with \eqref{adj} to yield $y^{\sharp\sharp} =
z^{\sharp\sharp} = 0$, forcing
\begin{align*}
e = S(y)^{-1}y^\sharp = S(z)^{-1}z^\sharp &&\text{(by
\eqref{yzz1})}  
\end{align*}
to be a primitive idempotent of $A$. The relation
\begin{align*}
e \times x = \,\,&S(y)^{-1}(y \times y^\sharp) + dS(z)^{-1}(z
\times z^\sharp) 
\end{align*}
combined with \eqref{pau}, \eqref{yzz2}, and \eqref{N0} gives:
\[
e \times x = -(y + dz) = T(x) (1 - e) - x.
\]
Lemma \ref{L-PEIRCEDEC}a shows that $x$ belongs to
$(A_K)_0(e)$, hence so does $\iota(x)$.

Now fix $v \in A_K$ such that $x = \iota(x) \times v$ and
write $v = ae + v_1 + v_0$ for $a \in K$ and  $v_i \in (A_K)_i(e)$ with $i =
0,1$. Since $\iota(x)$ is in $(A_K)_0(e)$, we may apply Lemma
\ref{L-PEIRCEDEC}a again to conclude $\iota(x) \times (ae) =
-a\iota(x) \in (A_K)_0(e)$. On the other hand, Lemma
\ref{L-PEIRCEDEC}c gives $\iota(x) \times v_0 =
S(\iota(x),v_0)e \in (A_K)_2(e)$.
Finally, using the circle product $a \circ b := \{ a, 1, b \} = U_{a,b}1$, we obtain
\[
\iota(x) \times v_1 = \iota(x) \circ v_1 -
T\big(\iota(x)\big)v_1 - T(v_1)\iota(x) + 
[T\big(\iota(x)\big)T(v_1) - T\big(\iota(x),v_1\big)]1
\]
by \eqref{bquatr} and \eqref{ads} linearized.  But, as the Peirce decomposition is orthogonal relative to the generic trace, this is just 
$\iota(x) \circ v_1$, which is in $(A_K)_1(e)$.
Decomposing $a = \alpha + \beta d$ with $\alpha, \beta \in k$ and
comparing Peirce components of $x = \iota(x) \times v$ relative to
$e$, a short computation yields
\begin{align*}
y + dz = \iota(x) \times v = -a\iota(x) = -\big(\alpha y + (\alpha
+ \delta\beta)z\big) - d(\beta y - \alpha z)\,,
\end{align*}
hence
\begin{align}
\label{COMPCOEF} (1 + \alpha)y = -(\alpha + \delta \beta)z \quad \text{and} \quad (1
- \alpha)z = -\beta y\,.
\end{align}

To complete the proof, it suffices to show $z = -2y$, because this
implies $0 = x^\sharp = (1 -2d)^2y^\sharp$, hence $2d = 1$, a
contradiction. Now, if $\beta = 0$, then $\alpha = 1$ and $2y =
-z$ by \eqref{COMPCOEF}. But if $\beta \neq 0$, then
\begin{align*}
\beta z^\sharp &= -\beta y \times z = (1 - \alpha)z \times z
&&\text{(by \eqref{yzz1} and \eqref{COMPCOEF})} \\
&= 2(1 - \alpha)z^\sharp\,,
\end{align*}
which yields $\beta = 2(1 - \alpha)$, hence $\beta z = 2(1 -
\alpha)z = -2\beta y$ by \eqref{COMPCOEF}, and again we end up
with $z = -2y$.
\end{proof}

\section{The case where $T(x) \ne 0$} \label{Tx1}

We now treat the other possible consequence of hypothesis (1).

\begin{prop} \label{Tx1.prop}
There is a nonzero element $x \in A_K$ satisfying
   \[
   x^\sharp = 0, \quad 
x \in \iota(x) \times A_K, \quad \text{and} 
\quad T(x) \ne 0
   \]
if and only if $k \times K$ is a subalgebra of $A$.   
\end{prop}

\begin{proof}
First suppose that $K$ is not a field, i.e., $K = k \times k$.  Then $A_K$ is identified with $A \times A$ and $\iota$ acts via the switch.  If $A_K$ contains an element $x = (x_1, x_2)$ as in the statement of the proposition, then one of the elements $x_1, x_2 \in A$ is singular, hence $A$ is reduced.  It follows that $k \times k \times k$ is a subalgebra.  Conversely, if $k \times k \times k$ is a subalgebra of $A$, then the standard basis vectors $e_1, e_2, e_3$ in $k^3$ form a complete orthogonal system of primitive idempotents in $A$.  Putting $x = (e_1, e_2)$, we find that $x^\sharp = 0$, $x$ has trace one, and 
\[
\iota(x) \times e_3 = (e_2 \times e_3, e_1 \times e_3) = x.
\]

We are left with the case where $K$ is a field.  Suppose first that $A_K$ contains an element $x$ as in the statement of the proposition.
Then putting $b := T(x) \in \Kx$, the element $e := b^{-1} x$ is a primitive idempotent in $A_K$ and
\begin{align*}
\iota(e) = \iota(b)^{-1}\iota(x) \in \iota(b)^{-1}(x \times A_K) = e
\times \big(\iota(b)^{-1}b \, A_K \big) \subseteq e \times A_K\,.
\end{align*}
But since $e \times e = 2e^\sharp = 0$ and $f = 1 - e$ is in $(A_K)_0(e)$, Lemma \ref{L-PEIRCEDEC}a implies that the map $x \mapsto e \times x$ kills $(A_K)_2(e) + (A_K)_1(e)$ and stabilizes $(A_K)_0(e)$.  So  $\iota(e)$ is in $(A_K)_2(e)$ and $e, \iota(e), c$ are orthogonal primitive idempotents in $A_K$ where $c
:= 1 - e - \iota(e)$.  The idempotent $c$, remaining fixed under $\iota$,
belongs to $A$. Thus
\begin{align*}
k \times K \ra A \quad \text{via}\ \quad (\alpha, a) \mapsto
\alpha c + ae + \iota(a)\iota(e)\,,
\end{align*}
is an embedding of Jordan algebras over $k$.

Conversely, suppose that $K$ is a field and $k \times K$ is a subalgebra of $A$.
By hypothesis, $ A\otimes K$ contains $K \times (K \otimes K) = K
\times K \times K$ as a subalgebra whose unit vectors
$u_1, u_2, u_3$ form a complete orthogonal system of primitive
idempotents in $A \otimes K$ such that $u_1$ is in $A$ and $\iota$ interchanges $u_2$ and $u_3$.
Thus $u_2 \in A \otimes K$ is a singular element with trace 1 that satisfies
\[
u_2 = u_3 \times u_1 \in \iota(u_2) \times A_K.\qedhere
\]
\end{proof}

\section{Sufficient condition for isomorphism} 

Let $A$ be a reduced Albert algebra over $k$ with
coordinate algebra $C$. By \ref{INVMOD2}, the 5-Pfister form
corresponding to $f_5(A)$ may be written as $N_C \otimes \gamma$
for some 2-Pfister bilinear form $\gamma$. Since reduced Albert
algebras are classified by their mod-2 invariants, we can write $A$ as $\hrm_3(C,\gamma)$. By the same
token, given octonion algebras $C$, $C'$ and 2-Pfister bilinear
forms $\gamma$, $\gamma'$ over $k$, the algebras
$\hrm_3(C,\gamma)$ and
$\hrm_3(C',\gamma')$ are isomorphic if and only
if $N_C \cong N_{C'}$ and $N_C \otimes \gamma \cong
N_{C'} \otimes \gamma'$. The goal of this section is
to prove a weak analogue of this statement for hermitian Albert
triples.

\begin{prop} \label{suff}
Let $K$ be a quadratic \'etale $k$-algebra, let $C$ be an octonion $k$-algebra, and let $\gamma$ and $\gamma'$ be 2-Pfister bilinear forms.  If $\gamma \ot N_K$ is isomorphic to $\gamma' \ot N_K$, then $\T(\hrm_3(C, \gamma), K)$ is isomorphic to $\T(\hrm_3(C, \gamma'), K)$.
\end{prop}

It is natural to wonder if a stronger result holds, namely if 
$N_C \ot N_K \ot \gamma$ is isomorphic to $N_C \ot N_K \ot \gamma'$, are the triples $\T(\hrm_3(C, \gamma), K)$ and $\T(\hrm_3(C, \gamma'), K)$ necessarily isomorphic?  (The tensor products in the question only make sense in characteristic different from 2.)  The answer is no, as Example \ref{aniso} below shows.

The proposition will follow easily (see the end of this
section) from an alternative construction of hermitian Albert
triples that we now describe. First we claim that the
general linear group $\GL_3(k)$ acts on the split Albert algebra
$A^d = \hrm_3(C^d, 1)$ via
\begin{align*}
\phi_g(j) := gjg^t &&\text{($g \in \GL_3(k), j \in A^d$)}
\end{align*}
in such a way that
\begin{align}
\label{SIM} N\big(\phi_g(j)\big) = (\det g)^2N(j)\,.
\end{align}
To see this, it suffices to show $gA^dg^t \subseteq A^d$ and
\eqref{SIM} for \emph{elementary} matrices $g \in \GL_3(k)$,
which follows easily by brute force using \eqref{NORM}. (This is
the argument given in \cite[\S5]{Jac:J3}.) Also, we have $\phi_g^\dag =
\phi_{g^{-t}}$ by the argument in \cite[p.~77]{Jac:J3}.

In particular, the map $\phi_g$, for $g \in \SL_3(k)$, is an
automorphism of the split hermitian Albert triple $\mathcal{T}^d$
via
   \[
   \phi_g \cdot (j_1, j_2) = (gj_1 g^t, g^{-t} j_2 g^{-1}).
   \]
This gives a map
   \begin{equation} \label{const.map}
   G_2 \times (\SL_3 \rtimes \Zm2) \ra \aut(\Td),
   \end{equation}
and a corresponding map 
   \begin{equation} \label{const.2}
   H^1(k, G_2) \times H^1(k, \SL_3 \rtimes \Zm2) \ra H^1(k, \aut(\Td)).
   \end{equation}
This last map takes an octonion $k$-algebra $C$, a quadratic \'etale $k$-algebra $K$, and a rank 3 $K/k$-hermitian form as inputs, and it gives a hermitian Jordan triple as output.

\begin{lem} \label{Ferr.lem}
For every quadratic \'etale $k$-algebra $K$, octonion $k$-algebra $C$, and 3-dimensional symmetric bilinear form $\gamma$, construction \eqref{const.2} sends the hermitian form deduced from $\gamma$ to $\T(\hrm_3(C, \gamma), K)$.
\end{lem}

\begin{proof}[Sketch of proof]
The proof is essentially the same as the proof of \cite[Prop.~3.2]{Ferr:E62}, so we omit it.  One observes that construction \eqref{const.2} sends the hermitian form to $\T(A, K)$, where $A$ is the $\gamma$-isotope of $\hrm_3(C, 1)$.  But $A$ is isomorphic to $\hrm_3(C, \gamma)$ by \cite[p.~61, Th. 14]{Jac:J} or \cite[p.~1077]{McC:NAS}.
\end{proof}

\begin{eg}
In \cite{Veld:unit1} and \cite{Veld:unit2}, Veldkamp considers groups of type 
$\dE$ constructed by \eqref{const.2} from an octonion $k$-algebra $C$ and 
a cocycle in $Z^1(K/k, \SL_3 \rtimes \Zm{2})$ whose value at $\iota$ is $(\qform{1, 1, -1}, 1)$.  This cocycle corresponds to the $K/k$-hermitian form deduced from
$\qform{1, 1, -1}$.  
By the lemma, such a group is isomorphic to $G(\hrm_3(C, \qform{1, 1, -1}), K)$, i.e., is one of the groups in the first three lines of Table \ref{iso.table}.
\end{eg}

\begin{proof}[Proof of Proposition \ref{suff}]
Under our hypotheses, the rank 3 $K/k$-hermitian forms deduced from the pure parts of $\gamma$ and $\gamma'$ are isometric.  Therefore, the hermitian Jordan triples constructed from them via \eqref{const.2} are isomorphic.  Now apply Lemma \ref{Ferr.lem}.
\end{proof}

\begin{borel}{\textsl{Proof that (3) implies (1)}} \label{iso.cor}
Suppose that (3) holds in the statement of Theorem \ref{MT}.  Then $A$ is isomorphic to $\hrm_3(C, \gamma)$ where $\gamma \ot N_K$ is hyperbolic.  By Prop.~\ref{suff}, $G(A, K)$ is isomorphic to $G(\hrm_3(C, \qform{1, -1, 1}), K)$, and (1) holds by Example \ref{nil.eg2}.\hfill$\qed$%
\end{borel}

\begin{borel}{Vista}
The center of $\SL_3$ is identified with the center of the split group $E_6 := \aut(\Td)$ via \eqref{const.map}, and this gives a map
   \[
   G_2 \times \aut(\PSL_3) \ra \aut(E_6).
   \]
This leads to a construction of groups of type $E_6$ (with possibly non-trivial Tits algebras), with inputs an octonion algebra, say $C$, and a central simple associative algebra of degree 3 with unitary involution fixing $k$, say $(B, \tau)$.  Tits's other construction of Lie algebras of type $E_6$ (corresponding to the  $E_6$ in the bottom row of the ``magic square'' \cite[p.~98]{Jac:ex}) uses identical inputs, and it is natural to guess that it produces the Lie algebra of the group arising from the construction above.  In any case, one can ask: \emph{What is the index of the group constructed from given $C$ and $(B, \tau)$?}  

We remark that some of the techniques in this paper can be adapted to attack this question.  For \S\ref{hjp}, hermitian Jordan triples should be replaced by a new type of algebraic structure in a manner completely analogous to the replacement of quadratic forms by algebras with orthogonal involution as in \cite[\S5]{KMRT}.  The description of the homogenous projective varieties in \S\ref{homo} can be translated directly to this new structure.  We leave the details to the interested reader.
\end{borel}

\section{The case where $A$ is split by $K$} 

\begin{borel*} \label{furthermore}
In the notation of the main theorem, we always have (4) implies (1) by Example \ref{nil.eg2}.  Conversely, suppose that $K$ splits $A$ and (3) holds.  In particular, $A$ is reduced and we may write $f_3(A)$ as $[K] \cdot \tau$ for some 2-Pfister bilinear form $\tau$, hence 
   \[
   f_5(A) = f_3(A) \cdot \gamma = \tau \cdot 0 = 0.
   \]
That is, $A$ has nonzero nilpotents.  This proves the final sentence of Theorem \ref{MT}.
\hfill$\qed$
\end{borel*}

\smallskip

We now provide an example to show that one really needs the hypothesis that $K$ splits $A$ in the implication that (3) implies (4).

\begin{eg}
Take $k = \Q(x, y, z, u, d)$ to be the rational function field in five variables over $\Q$, let $C$ be the octonion $k$-algebra corresponding to the 3-Pfister form $\pform{x, y, z}$, and put
   \[
   A := \hrm_3(C, \qform{d, 1, -u}) \quad \text{and} \quad K := k(\sqrt{d}).
   \]
Then $f_3(A) = (x) \cdot (y) \cdot (z)$ and $f_5(A) = f_3(A) \cdot \gamma$ for $\gamma = \pform{-d, u}$.  Since $(-d) \cdot (d) = 0$, we have $\gamma \cdot (d) = 0$, i.e., (3) holds.  On the other hand, the 5-Pfister $\pform{x, y, z, -d, u}$ is anisotropic by Springer's Theorem \cite[\S{VI.1}]{Lam}, so $f_5(A)$ is nonzero and $A$ contains  no nonzero nilpotents.
\end{eg}

\section{Final observations} 

We close this paper by applying our main theorem to give an easy criterion for when $G(A, K)$ is isotropic over special fields.

\begin{prop} \label{SAP}
Let $k$ be a SAP field of characteristic zero such that $\cd_2 k(\sqrt{-1}) \le 2$.  A group $G(A, K)$ is isotropic if and only if $A$ is reduced and $f_5(A) \cdot [K] = 0$.
\end{prop}

Every algebraic extension of $\Q$ (not necessarily of finite degree) and $\R((x))$ satisfy the hypothesis of the proposition.  See \cite[pp.~653--655]{BP:H} for a summary of basic properties of fields $k$ as in the proposition.

The restriction on the characteristic is harmless.  The prime characteristic analogue of the proposition is a corollary of the statement: \emph{If $H^5(k, 2)$ is zero and $A$ is reduced, then $G(A, K)$ is isotropic.}  To see this, note that the hypotheses imply that the subgroup $\aut(A)$ of $G(A, K)$ is itself isotropic.

\begin{proof}[Proof of the proposition]
If $G(A, K)$ is isotropic, then $A$ is reduced and $f_5(A) \cdot [K]$ is zero by the main theorem, so we assume that $A$ is reduced and $f_5(A) \cdot [K]$ is zero and prove the converse.  By the SAP property, there is some $\gamma_1 \in \kx$ that is positive at every ordering where $f_5(A)$ is zero and negative at every ordering where $f_5(A)$ is nonzero.  Put $\gamma := \pform{-1, \gamma_1}$.  Because $f_3(A)$ divides $f_5(A)$, it follows that $f_5(A)$ equals $f_3(A) \cdot \gamma$ and $\gamma \cdot [K]$ equals zero over every real-closure of $k$, which in turn implies those same equalities over $k$.  We conclude that $G(A, K)$ is $k$-isotropic by the main theorem.
\end{proof}

We observed in Prop.~\ref{suff} that 
one can change $\gamma$ somewhat without changing the isomorphism class of $G(\hrm_3(C, \gamma), K)$.  Motivated by Prop.~\ref{SAP}, one might hope that the main theorem still holds with (3) replaced by
\begin{enumerate}
\item[($3'$)] $A$ is reduced and $f_5(A) \cdot [K] = 0$.
\end{enumerate}
(We remark that the expression $f_5(A) \cdot [K]$ only makes sense when $\chr k \ne 2$.)
Clearly, (3) implies ($3'$).  We now give an explicit example where ($3'$) holds but (3) does not.  That is, such a replacement is not possible.   We thank Detlev Hoffmann for showing us this example.

\begin{eg} \label{aniso}
Fix $k_0$ to be the purely transcendental extension $\Q(x, y, z, u, v, d)$ of the rationals (say).  Let $k$ denote the function field over $k_0$ of the 6-Pfister form $\pform{x, y, z, u, v, d}$.  Define a $k$-group $G(A, K)$ via
  \[
  A := \hrm_3(C, \qform{-u, -v, uv}) \quad \text{and} \quad K := k(\sqrt{d}),
  \]
where $C$ is the octonion $k$-algebra with norm $\pform{x, y, z}$.  
Clearly, ($3'$) holds over $k$.

On the other hand, take $E$ to be the function field of $q := \pform{x, y, z}' + \qform{d}$ over $k$, where $\pform{x, y,z}'$ denotes the pure part of the 3-Pfister.  Clearly, $f_3(A) = (x) \cdot (y) \cdot (z)$ is killed by $EK$, that is, $[K]$ divides $f_3(A)$ over $E$.  
We now argue that $f_5(A)$ is not zero over $E$, which will show that (3) fails over $E$ hence also over $k$.  Note that $f_5(A) = (x) \cdot (y) \cdot (z) \cdot (u) \cdot (v)$ is nonzero over $k_0$ and remains nonzero over $k$ by Hoffmann's Theorem \cite[X.4.34]{Lam}. 

For sake of contradiction, suppose that $f_5(A)$ is killed by $E$.  Then by Cassels-Pfister \cite[X.4.8]{Lam}, for $a \in \kx$ represented by $q$ and $b \in \kx$ represented by $f_5(A)$, we have that $\qform{ab} q$ is a subform of $f_5(A)$.  Taking $a = b = -x$, we find that $q$ is a subform of $f_5(A)$ over $k$.  Computing in the Witt ring of $k$, it follows that 
   \[
   \qform{1, -d} + \pform{x, y, z} \qform{-u, -v, uv}
   \]
is $k$-isotropic.  But this form is $k_0$-anisotropic and remains anisotropic over $k$ because its dimension is $26 < 32 < 64 = \dim \pform{x, y, z, u, v, d}$, again using Hoffmann's Theorem.  This is a contradiction, which shows that (3) does not hold.
\end{eg}

{\small{\subsection*{Acknowledgements} The authors are indebted to Detlev Hoffmann, Ottmar Loos, Erhard Neher, and Michel Racine for useful discussions on the subject of this paper.}}


\providecommand{\bysame}{\leavevmode\hbox to3em{\hrulefill}\thinspace}
\providecommand{\MR}{\relax\ifhmode\unskip\space\fi MR }
\providecommand{\MRhref}[2]{%
  \href{http://www.ams.org/mathscinet-getitem?mr=#1}{#2}
}
\providecommand{\href}[2]{#2}

\end{document}